%% file: hal.tex
\documentclass[11pt]{article}
\pdfoutput=1 
\usepackage[utf8x]{inputenc}
\usepackage[english]{babel}
\usepackage{graphicx}           
\usepackage{url}
\usepackage{eurosym}
\usepackage{makeidx}
\usepackage[plain]{algorithm}
\usepackage{algorithmic} 
\usepackage{color} 
\usepackage[margin=3.3cm]{geometry}
\usepackage{setspace}
\usepackage{datetime}
\usepackage{pifont}
\usepackage{tikz}
\usepackage[numbers]{natbib}
\usepackage{mathtools}
\usepackage{amsmath,amsfonts,amssymb,latexsym,stmaryrd}

\usepackage{import}
\usepackage{transparent}

\begin{document}

\newcommand{\Tmax}{T^{\textrm{\tiny max}}}
\newcommand*\circled[1]{%
  \tikz[baseline=(C.base)]\node[draw,circle,inner sep=0.5pt](C) {#1};\!
}
\newcommand{\circledi}{{\scriptsize $\protect\circled{\it i\/}$}}
\newcommand{\circledip}{{\scriptsize $\protect\circled{\it i\/}\;'$}}

\newcommand{\tto}{t^\textrm{\rm\tiny o}}
\newcommand{\xio}{\xi^\textrm{\rm\tiny o}}

\newcommand{\QP}{Q^{\textrm{\tiny P}}}
\newcommand{\QB}{Q^{\textrm{\tiny B}}}
\newcommand{\QE}{Q^{\textrm{\tiny E}}}
\newcommand{\qP}{q^{\textrm{\tiny P}}}
\newcommand{\qB}{q^{\textrm{\tiny B}}}

\newtheorem{theorem}      {Theorem}[section]
\newtheorem{lemma}        [theorem]{Lemma}
\newtheorem{remark}       [theorem]{Remark}

\newcommand{\E}        {\mathbb E}
\newcommand{\N}        {\mathbb N}
\renewcommand{\P}      {\mathbb P}
\newcommand{\R}      {\mathbb R}
\newcommand{\CC}   {{\mathcal C}}
\newcommand{\LL}   {{\mathcal L}}
\newcommand{\NN}   {{\mathcal N}}
\newcommand{\AAA}  {{\mathcal A}}

\newcommand{\eqdef}     {\stackrel{{\textrm{\rm\tiny def}}}{=}}
\newcommand{\simiid}     {\stackrel{\textrm{iid}}{\sim}}
\newcommand{\equiva}    {\displaystyle\mathop{\simeq}}
\newcommand{\demi}{{{\textstyle\frac{1}{2}}}}
\newcommand{\indic}{{\mathrm\mathbf1}}
\newcommand{\rmd}   {{{\textrm{\upshape d}}}}
\newcommand{\proof}        {\paragraph{Proof}}


\title{Estimation of the parameters\\ of a stochastic logistic growth model}
\author{Fabien Campillo\thanks{\protect\url{Fabien.Campillo@inria.fr} --- 
           Project--Team MODEMIC, INRIA/INRA, UMR MISTEA, b\^at. 29, 2
           place Viala, 34060 Montpellier cedex 06, France.}  
   \and Marc Joannides\thanks{\protect\url{marc.joannides@univ-montp2.fr} ---
   Universit\'e Montpellier 2 / I3M, case courrier 51,
        place Eug\`ene Bataillon, 34095 Montpellier cedex 5; this
        author is associate researcher for Project--Team MODEMIC,
        INRIA/INRA, UMR MISTEA.} 
   \and Irène Larramendy-Valverde\thanks{\protect\url{irene.larramendy-valverde@univ-montp2.fr} ---
   Universit\'e Montpellier 2 / I3M, case courrier 51,
        place Eug\`ene Bataillon, 34095 Montpellier cedex 5.}
      }
\date{}

\maketitle

\begin{abstract}
We consider a stochastic logistic growth model involving both birth
and death rates in the drift and diffusion coefficients for which
extinction eventually occurs almost surely. The associated complete 
Fokker--Planck equation describing the law of the process is
established and studied. We then use its solution to build a
likelihood function for the unknown model parameters, when discretely
sampled data is available. The existing estimation methods need
adaptation in order to deal with the extinction problem. We propose
such adaptations, based on the particular form of the Fokker--Planck 
equation, and we evaluate their performances with numerical
simulations. In the same time, we explore the identifiability of the
parameters which is a crucial problem for the corresponding
deterministic (noise free) model. 
\paragraph{Keywords and phrases:} 
Logistic model, diffusion processes, extinction, Fokker--Planck
equation, estimation, Monte Carlo.
\end{abstract}

\section{Introduction}
Most of the growth models in population dynamics and ecology are based
on ordinary differential equations (ODE). Among these models, the
logistic population growth model was first introduced
by~\citet{verhulst1838a} to take into account crowding effect, by 
damping the {\it per capita growth rate} in the Malthusian growth
model; this model reads:
\begin{equation}
  \label{eq.verhulst} 
  \dot x(t) = r\, x(t) \, \left ( 1 - \frac{x(t)}{K}\right ) \,,\qquad
  x(0) = x_0 > 0 \,,\ t\geq 0
\end{equation} 
where $x(t)\geq 0$ is the density of some population, $r>0$ the growth
rate and $K>x_0$ the carrying capacity of the environment. This
formulation leaves aside the stochastic features resulting from
diversity in the population or from random fluctuations of the
environment.

Parameters in \eqref{eq.verhulst} are usually identified using least
squares methods, but this approach obscures two issues. On the one
hand the model \eqref{eq.verhulst} cannot account for a possible
extinction of the population and therefore cannot benefit from the
information in a data set with extinction. On the other hand only the
growth rate $r=\lambda-\mu$, which is the difference between the birth
rate $\lambda$ and death rate $\mu$, can be identified in this context
and~\eqref{eq.verhulst} cannot provide information on each of these
rates separately.

\medskip

Stochastic counterparts of \eqref{eq.verhulst}, usually expressed as
stochastic differential equations (SDE), may overcome these two issues.
The stochastic model that we consider explicitly handles the question
of extinction; it also makes the information contained in the
demographic noise available, leading to a rough approximation of
$\lambda+\mu$. Stochastic logistics models can be obtained by adding a
random ad hoc perturbation term in~\eqref{eq.verhulst}. A more natural
way is to consider a diffusion approximation of a birth and death
process that features a logistic mechanism, see Appendix
\ref{appendix.sde}.  For both approaches, there will obviously be many
stochastic models derived from or leading to the same deterministic
model \eqref{eq.verhulst}, but with different qualitative behaviors,
see~\citet{schurz2007a}.

In this paper, we will consider the stochastic logistic model given by
the following SDE:
\begin{align}
  \label{eq.sde}
  \rmd X_{t} &= (\lambda-\mu-\alpha\,X_{t})\,X_{t}\,\rmd t + \rho\,
  \sqrt{(\lambda+\mu+\alpha\,X_{t})\,X_{t}}\,\rmd B_{t}
\end{align}
where $\lambda>0$ is the birth rate, $\mu>0$ the death rate,
$\alpha>0$ the logistic coefficient and $\rho>0$ the noise intensity
which relates to the order of magnitude of the underlying population
(see Appendix \ref{appendix.sde}); $B_{t}$ is a standard Brownian
motion; the law of the initial condition $X_0$ is supported by
$\R_{+}$; $B_{t}$ and $X_{0}$ are supposed independent.

The objective of this paper is to study the estimation problem of
model~\eqref{eq.sde} for the unknown  parameter 
$\theta\in\Theta=(0,\infty)^p$, based on a discrete sample of one
trajectory. The parameter  $\theta$ may include some or all parameters
$(\lambda,\mu,\alpha,\rho)$ and may also appear in the initial
distribution law. Hence~\eqref{eq.sde} can be rewritten: 
\begin{align}
\label{eq.sde2}
  \rmd X_{t}
  &=
  b^\theta(X_{t})\,\rmd t+\sigma^\theta(X_{t})\,\rmd B_{t}\,,
  \quad
  0\leq t\leq T\,,\quad X_{0}\sim \pi_0^\theta(\rmd x)
\end{align}
where $b^\theta(x) = (\lambda - \mu - \alpha \, x)\,x$ and $\sigma^\theta (x) = \rho\,
\sqrt{(\lambda + \mu + \alpha \, x )\,x}$
are  the drift and diffusion coefficients respectively; $\pi_0^\theta$
is the initial distribution law.  We also define $a^\theta (x) \eqdef
[\sigma^\theta (x)]^2$.

Due to the Markov nature of the
process given by \eqref{eq.sde2}, the distribution law of the data can
be expressed as a product of the transition kernel between successive
instants of observation. The latter is known to be given, in a weak
form, by the Kolmogorov Forward Equation. For diffusion process that
never becomes extinct, this equation reduces to the Fokker-Planck
Equation for the transition density. An originality of this work
lies in the fact that the solution of this Kolmogorov equation fails
to have a density with respect to the Lebesgue measure on~$\R_+$. We
investigate the complete form of the Fokker-Planck Equation
(\citet{feller1952a}) that gives the evolution of the transition
kernel of the diffusion process 
$\{X_{t} \}_{0\leq t\leq   T}$ in Section~\ref{sec.statistical.model}:
in Section~\ref{sec.statistical.model.exit} we prove  that  $x=0$ is
an exit boundary point according to Feller terminology;   in
Section~\ref{sec.statistical.model.FP} we establish the Fokker-Planck
Equation and finally the likelihood function is detailed in
Section~\ref{sec.statistical.model.LF}. The transition kernel and the
likelihood function derived from it in
Section~\ref{sec.statistical.model} cannot be computed explicitly. In
Section~\ref{sec.approximation} we propose to adapt existing numerical 
approximation procedures to take extinction into account: in
Section~\ref{sec.approximation.fpe} we develop specific finite
difference schemes in order to approximate the solution Fokker--Planck
equation; 
in Section~\ref{sec.approximation.MC} we propose appropriate Monte
Carlo approximations. The properties of the model and the
approximations methods performances are numerically investigated in
Section~\ref{sec.numerical.experiments}. Appendices are dedicated to
the development of the logistic SDE, to the  existence and uniqueness
of solutions of the SDE, and to the algorithmic description of the
Monte Carlo methods considered.

\section{Statistical model}
\label{sec.statistical.model}

In Appendix \ref{appendix.sde} we prove that (\ref{eq.sde}) admits a
unique solution. In the present section we describe the nature of the
boundary point $0$ and we establish the Fokker-Planck equation that
gives the evolution of the transition kernel of the diffusion
process~\eqref{eq.sde2}. This Fokker-Planck equation explicitly
handles the 
probability of extinction. The statistical model and the likelihood
function are derived at the end of this section. 
For notational simplicity, we drop the reference of the
parameter $\theta$ in the two next subsections.

\subsection{Extinction time}
\label{sec.statistical.model.exit}

For $y \geq 0$, let:
\[
 \tau_y \eqdef \inf\{t \geq 0\,; \,  X_t = y \}
\]
As a
by--product of the proof of the  existence and uniqueness of solutions
of \ref{eq.sde2} given in Appendix~\ref{appendix.sde}, we find that
the process remains 
in the interval $[0,+\infty)$. We also show that $X_t = 0$ for
$t\geq \tau_0$, but whether the boundary~0 could be reach in finite
time or not  is still to be determined. A complete description of the
possible behavior at the boundary points has been established
by~\citet{feller1952a}. A detailed review of these results can be 
found in Chapter 15 of~\citet{karlin1981a}. The
following lemma states that~0 is an \emph{exit} boundary point
according to Feller terminology: it is reached in an almost surely
finite time and no interior point in $(0;+\infty)$ can be reached
starting from~0. 

\begin{lemma} 
\label{lemma.extinction}
Extinction occurs almost surely in finite time, that is for all $x
\geq 0$, $\P_x(\tau_0 < \infty) = 1$ where $\P_{x}$ is the probability measure such that $X_{0}=x$. 
\end{lemma} 

\proof
For $0 < x_l < x < x_r$, we have 
\begin{displaymath}
  \P_x (\tau_{x_l} < \tau_{x_r}) 
  = 
  \frac{S(x_r) - S(x)}{S(x_r) - S(x_l)}
\end{displaymath} 
where $S$ is the \emph{scale function}, see
e.g.~\citet{klebaner2005b}, defined by 
\begin{displaymath}
  S(x) 
  \eqdef 
  \int_\eta ^x 
    \exp\left\{ 
        -\int_\eta ^y \frac{2\, b(z)}{a(z)}\rmd z 
    \right \} \rmd y \,.
\end{displaymath} 
The choice of the lower bound $\eta$ in the integrals will appear to
be unimportant and could be chosen arbitrarily within $(x_l,x_r)$
since the relevant expressions involve only differences of the
function~$S$. A straightforward computation gives for this particular
case 
\begin{displaymath}
   S(x) = C_\eta\, \int_\eta ^x \,s(y)\, \rmd y
\end{displaymath} 
where
\begin{displaymath}
  s(y) 
  \eqdef 
  \left( 
    \mathrm{e}^y \, (\lambda + \mu + \alpha \, y)
          ^{-\frac{2\,\lambda }{\alpha}} 
  \right)^\frac{2}{\rho ^2}  
\end{displaymath}
and $C_\eta$ is a constant depending on $\eta$ only, so that
\begin{displaymath} 
  \P_x (\tau_{x_l} < \tau_{x_r}) 
  = 
  \frac{
   \int_\eta ^{x_r} s(y)\, \rmd  y - 
   \int_\eta ^x s(y)\, \rmd  y }
  { \int_\eta ^{x_r} s(y)\, \rmd  y - 
   \int_\eta ^{x_l} s(y)\, \rmd  y}
   = 
   1 -  \frac{ \int_{x_l} ^x s(y)\, \rmd  y} 
             { \int_{x_l} ^{x_r} s(y)\, \rmd y} \,.
\end{displaymath} 
Taking the limit as $x_l \downarrow 0$ yields
\begin{displaymath} 
  \P_x (\tau_0 < \tau_{x_r}) 
  = 
  1 - \frac{ \int_0 ^x s(y)\, \rmd  y} 
           { \int_0 ^{x_r} s(y)\, \rmd y} 
\end{displaymath} 
where both integrals are finite since $s$ is continuous on the compact
$[0,x_r]$. For the same reason, we  have $\lim_{x_r\uparrow \infty}
\int_0 ^{x_r} s(y)\, \rmd y = \infty$. Note also that we
already have from Appendix \ref{appendix.sde.existence.uniqueness},
$\lim_{x_r \uparrow \infty} \tau_{x_r} = \infty$, a.s. since explosion
does not occur. The probability of ultimate extinction is then 
\begin{displaymath}
   \P_x (\tau_0 < \infty) 
   = 
   \lim_{{x_r}\uparrow \infty} \P_x (\tau_0 < \tau_{x_r}) = 1 \,.
\end{displaymath}
\hfill $\square$ 

\subsection{Fokker-Planck equation}
\label{sec.statistical.model.FP}

We denote by\footnote{Let $K$ and $K'$ be two transition kernels on
$\R_+$. Throughout this paper, we use the following notations:
\begin{itemize}
\item left action on test function: 
$Kf(x) \eqdef  \int_{\R_+} f(y)\, K(\rmd y\,\vert\, x)
$,
\item right action on measure:
$(\nu K)(\rmd y) \eqdef \int_{\R_+} \nu(\rmd x) \, K(\rmd y\,\vert\, x)
$,
\end{itemize} }:
\begin{displaymath} 
  Q_t(\rmd y\,\vert\, x) \eqdef \P(X_{s+t} \in \rmd y | X_s = x)
\end{displaymath} 
the transition  kernel of the Markov process $\{X_t\}_{t\geq 0}$  and
by $\pi_t(\rmd y) =  (\pi_0 Q_t)(\rmd y)$ the distribution of
$X_t$. We note that $Q_t(\rmd y\,\vert\, x) $ is not absolutely
continuous with respect to the Lebesgue measure on $\R_+$, because it
gives positive probability to the boundary point~$0$. The Lebesgue
decomposition of $Q_t(\cdot \,\vert\,x)$ into absolutely continuous
and singular parts reads  
\begin{equation} 
\label{eq.law.xt}
Q_t(\rmd y\,\vert\, x)  =  E_t(x)\, \delta_0(\rmd y) +
p_t(y\,\vert\,x) \,\rmd y \ . 
\end{equation} 
The transition kernel $Q_t(\rmd y\,\vert\, x)$ is a probability
measure for any $x \geq 0$, so that the \emph{extinction probability}
starting from~$x$ is  
\begin{displaymath}
E_t(x) = 1 - \int_0^\infty p_t(y\,\vert\,x) \,\rmd y \ . 
\end{displaymath} 
The transition kernel $Q_t(\rmd y\,\vert\, x)$ is absolutely
continuous with respect to the reference measure  on $\R_{+}$
\[
  m(\rmd y) \eqdef \delta_0(\rmd y) + \rmd y
\]
with density
\begin{subequations}
\label{eq.density} 
\begin{equation}
\label{eq.transition.density} 
  q_t (y\,\vert\,x) 
  \eqdef 
  \begin {dcases}
    E_t(x) \,,& \text{ if } y = 0\,,\\
    p_t (y\,\vert\,x) \,,&  \text{ otherwise.} 
  \end {dcases} 
\end{equation} 
We suppose also that the initial distribution $\pi_{0}$ is absolutely
continuous with respect to the reference measure $m(\rmd y) $, and we
let:
\begin{align}
\label{eq.initial.density} 
  q_{0}(y) 
  &\eqdef 
  \frac{ \pi_0 (\rmd y)}{  m(\rmd y)} 
  =
  \begin {dcases}
    E_0 \,,& \text{ if } y = 0\,,\\
    p_0 (y) \,,&  \text{ otherwise.} 
  \end {dcases} 
\end{align}
\end{subequations}

\medskip

We now establish the evolution equations for $E_t(x)$ and  $p_t
(y\,\vert\,x)$, for any $x>0$ fixed. Note that for $x=0$, $p_t
(y,\vert\,x)=0$  and $E_t(x)=1$ for all $t\geq 0$ and all $y \in
\R_+$. The Kolmogorov forward equation describes the evolution of
$Q_t$ in a weak sense:  
\begin{equation} 
\label{eq.kf.diff}
  \frac{\rmd }{\rmd t} Q_tf
  =  Q_t(\AAA f)
  \,,\quad \forall f\in \CC^\infty_{K}(\R_{+})
\end{equation} 
where $\AAA$ is the infinitesimal generator defined by:
\begin{equation}
\label{eq.gen:X} 
   \AAA f(x) 
   \eqdef
    b(x)\,f'(x) + \demi  \,a(x)\, f''(x)
\end{equation} 
and $\CC^\infty_{K}(\R_{+})$ is the set of functions differentiable for 
all degrees of differentiation and with compact support included in
$[0,+\infty)$. Using decomposition~\eqref{eq.law.xt}, 
\begin{align*} 
  Q_t(\AAA f)(x)
  &= 
  \AAA f(0) 
  + \int_0^\infty   b(y)\, p_t(y\,\vert\,x)\, f'(y)\,\rmd y
  + \frac12 \, \int_0^\infty  a(y) \, p_t(y\,\vert\,x)\,
    f''(y)\,\rmd y \\
  &=  
  \int_0^\infty  b(y)\, p_t(y\,\vert\,x)\,f'(y)\,\rmd y
  + \frac12 \, \int_0^\infty  a(y) \, p_t(y\,\vert\,x)\, 
           f''(y)\,\rmd y \,.
\end{align*} 
Note that $\AAA f(0)=0$ since both the drift and diffusion terms vanish
at~$0$.  A first integration by parts gives
\begin{align*}
  \int_0^\infty  b(y)\, p_t(y\,\vert\,x)\,f'(y)\,\rmd y
  & = 
  \big[b(y)\, p_t(y\,\vert\,x)\,f(y) \big]_0^\infty 
  -
  \int_0^\infty \, \frac{\partial ( b\, p_t(\cdot\,\vert\,x))
  }{\partial y}(y)\, f(y)\,\rmd y  
\\
  & = 
  -\int_0^\infty \,
  \frac{\partial ( b(y)\, p_t(y\,\vert\,x)) }{\partial y}\, f(y)\,\rmd
  y \,, 
\end{align*} 
and similarly
\begin{displaymath} 
  \int_0^\infty  a(y) \, p_t(y\,\vert\,x)\, f''(y)\,\rmd y \
  =  
  -\int_0^\infty \frac{\partial ( a(y)\, p_t(y\,\vert\,x)) }{\partial y}
        \,f'(y)\,\rmd y 
\end{displaymath} 
by the same property. In the above integrals, the non-integral terms
vanish at~$\infty$ because~$f\in \CC^\infty_{K}(\R_{+})$, but
they vanish at~$0$ because~$b(0) = a(0) = 0$. A second
integration by parts gives
\begin{align*} 
  &
  -\int_0^\infty \,
     \frac{\partial ( a(y)\, p_t(y\,\vert\,x)) }{\partial y}\,f'(y)\,\rmd y 
\\
  &\qquad\qquad = 
  - \left [  \frac{\partial ( a(y)\, p_t(y\,\vert\,x)) }{\partial y}\,
    f(y)  \right ]_0^\infty 
  + \int_0^\infty \,
  \frac{\partial^2 ( a(y)\, p_t(y\,\vert\,x)) }{\partial y^2}\, f(y)\,\rmd y 
\\
  &\qquad\qquad = 
  \left.\frac{\partial ( a(y)\, p_t(y\,\vert\,x)) }{\partial y}\right|_{y=0}\,f(0) 
 + \int_0^\infty \,
  \frac{\partial^2 ( a(y)\, p_t(y\,\vert\,x)) }{\partial y^2}\,
  f(y)\,\rmd y\,.
\end{align*} 
We define  $\AAA^\ast$ is the formal adjoint operator of~$\AAA$ acting on the
``forward'' space variable~$y$ only by
\begin{displaymath} 
 \AAA^\ast p_t(y\,\vert\,x) 
 = 
 - \frac{\partial [b(y)\,p_t(y\,\vert\,x)]}{\partial y }
 + \frac12 \, \frac{\partial^2 [a(y)\, p_t(y\,\vert\,x)]}{\partial y^2 }
 \,,
\end{displaymath} 
and we finally have the decomposition
\begin{equation}
\label{eq.Qt.af} 
 Q_t(\AAA f)(x)
 = 
 \frac12 \,  
    \left.\frac{\partial( a(y)\, p_t(y\,\vert\,x)) }{\partial y}\right|_{y=0}\,
    f(0)
+  \int_0^\infty \, \AAA^\ast p_t(y\,\vert\,x)\, f(y)\,\rmd y \,.
\end{equation} 
In view of~\eqref{eq.kf.diff},  the
first term of this decomposition has a nice interpretation: it is the
rate of increase of the extinction probability at time~$t$, expressed
as a probability flux through the boundary~$0$ (up to a minus
sign). Indeed, considering test functions $f_\epsilon$, such that
$f_\epsilon(0) = 1$, $f_\epsilon (y) = 0$ for $y\geq \epsilon$ and
with first two derivatives vanishing at~$0$, we get
\begin{displaymath}
  \frac{\rmd}{\rmd t} 
  \left[ E_t(x) + \int_0^\infty p_t(y\,\vert\,x)\, f_\epsilon(y)\; \rmd y \right] 
  = 
  \frac12 \, 
  \left.\frac{\partial ( a(y)\, p_t(y\,\vert\,x)) }{\partial y}\right|_{y=0}
  + 
  \int_0^\infty \, \AAA^\ast p_t(y\,\vert\,x)\, f_\epsilon (y)\,\rmd y \,.
\end{displaymath} 
The integrals vanish as $\epsilon \downarrow 0$ so that we obtain the
differential equation satisfied by $E_t(x)$:
\begin{subequations}
\label{eq.fpe}
\begin{equation}
\label{eq.fpe.E}
  \frac{\rmd}{\rmd t} E_t(x) 
  = 
  \frac12 \, 
  \left.
    \frac{\partial\big(a(y)\, p_t(y\,\vert\,x)\big) }
         {\partial y}
  \right|_{y=0}\, p_t(0\,\vert\,x) 
  \,,\quad E_0(x)=0\,.
\end{equation} 
On the other hand, the Fokker--Planck equation for the absolutely
continuous part ~$p_t(\cdot\,\vert\,x)$ is obtained by considering test functions
vanishing at~$0$ in~\eqref{eq.Qt.af}:
\begin{equation} 
\label{eq.fpe.p}
  \frac{\partial p_t(y\,\vert\,x)}{\partial t} =  \AAA^\ast p_t(y\,\vert\,x) \,, 
  \quad
   \lim_{t\downarrow 0}p_t(y\,\vert\,x) \,\rmd y=  \delta_{x}(\rmd y)
\end{equation} 
\end{subequations}
which is a PDE in a classical sense describing the evolution of the
process \emph{before extinction}. It follows that $y\mapsto p_t(y\,\vert\,x)$ is the
density of a \emph{defective distribution}. This equation has been
extensively studied by~\citet{feller1952a}. A notable result of
the latter work is that no boundary condition at~$0$ is required
for~\eqref{eq.fpe.p} to have a unique solution in~$L^1$. In
chapter~5 and~6 of~\citet{schuss2010a}, the multidimensional
case is investigated. In~\citet{campillo2013a}, the
authors establish an equation similar to~\eqref{eq.fpe} for a
two--dimensional model of a bioreactor, where extinction concerns
only one of the components.  
\begin{remark} 
According to Lemma~\ref{lemma.extinction}, $E_t(x)$ increases to~$1$, so
that $Q_t(\cdot\,\vert\,x)$ will eventually degenerate to the Dirac mass
at~$0$. We note that this convergence may be slow, i.e. that the
contribution of the Dirac mass in~\eqref{eq.law.xt} may not be
significant for the time scale at which the system is observed. This
phenomenon is investigated in~\citet{grasman1999asymptotic}.
\end{remark}

\subsection{Likelihood function}
\label{sec.statistical.model.LF}

We denote by $\P^{\theta}$ the underlying distribution of the process
$\{X_{t} \}_{0\leq t\leq   T}$.  Observations from the
SDE \eqref{eq.sde2} are available under the form: 
\[
  \xio_{k} = X_{\tto_{k}}\,,\quad k=0,\dots,M
\]
where, for sake of simplicity, the observation instants are equally
spaced, i.e. $\tto_{k}=k\,\Delta$ with $\Delta\eqdef T/M$.  By the
Markov property, the distribution of the measurements vector is
\begin{displaymath} 
  \P^\theta(\xio_0 \in \rmd \xi_0,\dots,\xio_M \in \rmd \xi_M) 
  = 
  \pi_0^\theta(\rmd \xi_0)\, \prod_{k=0}^{M-1} Q^\theta_\Delta(\rmd
  \xi_{k+1}\,\vert\,\xi_k) \,.
\end{displaymath} 
From the previous section, we know that $Q^\theta _\Delta(\rmd
y\,\vert\,x)$ 
has decomposition
\begin{align*} 
Q^\theta_\Delta(\rmd y\,\vert\,x) 
 = 
 E^\theta_\Delta(x)\, \delta_0(\rmd y) +  
 p^\theta_\Delta(y\,\vert\,x)\, \rmd y
\end{align*} 
where $p^\theta_t(x,y)$ and $E^\theta_t(x)$ solve the Fokker--Planck
equation~\eqref{eq.fpe}. 

The transition kernel $Q^\theta _\Delta(\rmd y\,\vert\,x)$ and the initial
distribution $\pi^\theta_{0}(\rmd y)$ are absolutely continuous with
respect to the reference measure  $m(\rmd y) = \delta_0(\rmd y) + \rmd
y$  with density $q^\theta _\Delta(y\,\vert\,x)$ and the initial distribution
$q^\theta_{0}(y)$ given by \eqref{eq.density}. Our statistical model
is therefore dominated by the product measure $m(\rmd \xi_0)\dots
m(\rmd \xi_M)$, and a likelihood function is given by  
\begin{align}
\label{eq.likelihood}
  \LL(\theta)
  &=
  q^\theta_{0}(\xi_{0})\,
  \prod_{k=0}^{M-1} q^\theta_{\Delta}(\xi_{k+1}\,\vert\,\xi_{k})\,.
\end{align}
The main object of interest is $q^\theta_\Delta(y\,\vert\,x)$ for given~$x$,
which is the solution of the set of differential
equations~\eqref{eq.fpe} for $t\in[0,\Delta]$. There is no explicit
solution available for this equation; we will therefore rely on
different types of approximations, either numerical or
analytical. This is the subject of next section.
%
%
%
%
%

\section{Transition kernel approximations}
\label{sec.approximation}

Notice that the density searched for consists
of two distinctive parts. The continuous component~$p_t$ can be
approximated independently whereas the discrete component~$E_t$
strongly depends on~$p_t$. This suggests that we must first design an
approximation to~$p_t$ from which the approximation of~$E_t$ can be
deduced. In all cases, any acceptable approximation should be a
probability density. Throughout this section, we will drop the
index~$\theta$ which is of no use for the description of the methods.

\subsection{Finite difference approximation of the Fokker-Planck equation}
\label{sec.approximation.fpe}

Numerical approximations of~\eqref{eq.fpe.p} can be obtained by
classical methods of numerical analysis of PDEs, paying attention to
the specific features of our model. Indeed, any appropriate
discretization scheme should correctly handle the degeneracy
(vanishing diffusion) at~$0$. Also the approximated solution should
remain non negative and integrate to at most~$1$, since it approaches a
defective probability distribution. Finally, the mass default must
be a consistent approximation of~\eqref{eq.fpe.E}. The approach
presented in~\citet{kushner1992a} seems natural in
this context, because it allows a straightforward interpretation of
the discretized operator in terms of generator of a Markov
process. See~\citet{campillo2013a} for such a
discretization method applied to a two--dimensional model. 

\subsubsection*{Space discretization}

\begin{figure}[h]
  \centering
\includegraphics[width=0.58\linewidth]{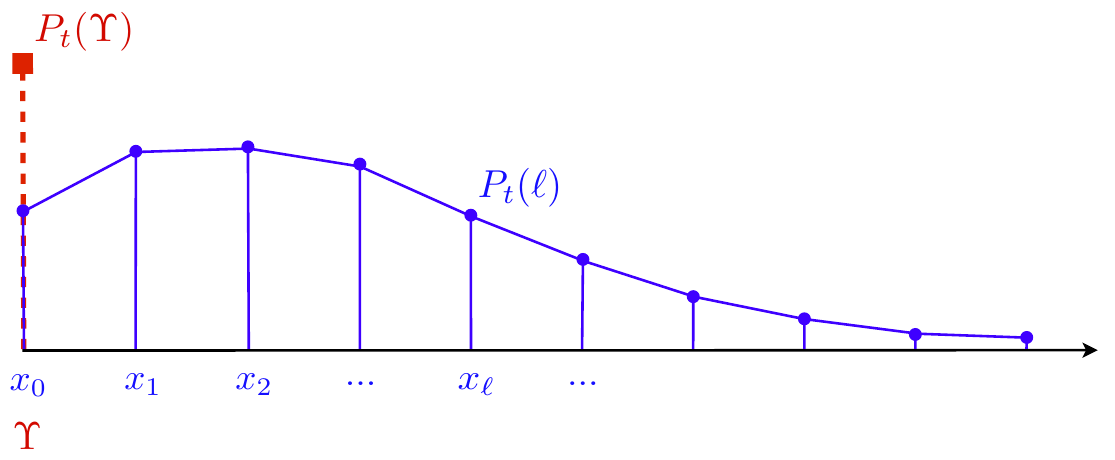}
\caption{Discretization of the state space as a regular finite
  grid. Value~$y=0$ is either the node~$x_0$ at which the value
  of~$p_t$ is evaluated and the cemetery point~$\Upsilon$.}
\label{fig.space}
\end{figure}

We discretize the space as a regular grid:
\[
    x_\ell = \ell\, h\,,\quad \ell=0,\dots,L
\]
for $h$ and $L$ given. Note that this grid is finite so
that it does not cover the whole support of $p_t(\cdot\,\vert\,x)$. In numerical
experiments, the range of the grid will have to be large enough so
that any artificial boundary condition imposed at $x_L$ will cause
limited harm. More importantly, the boundary point~$0$ has a twofold
status; as the node~$x_0$ of the grid, it enters the computation
of the continuous component~$p_t(0\,\vert\,x)$ and as an absorbing state, it
carries the extinction probability~$E_t$. It is thus legitimate to
introduce an additional~\emph{cemetery point}~$\Upsilon$ at
location~$0$, see Figure~\ref{fig.space}. Indeed, such a decomposition
of the point~$0$ gives the expected smoothness of~$p$ at the boundary,
depicted on Figure~\ref{fig.space} and observed on Figure~\ref{fig.kernel}.

%
We now derive the finite difference approximation of the continuous
part~$p_t$, returning to the weak formulation.  For suitable test
function $\phi$, 
\begin{align*} 
\int_0 ^\infty p_t(y\,\vert\,x)\, \AAA \phi(y)\, \rmd y 
&\simeq \frac{h}{2} \, p_t(0\,\vert\,x)\, \AAA \phi(0) +   
  h\, \sum_{\ell = 1}^{L -  1} p_t(x_\ell \,\vert\,x)\, \AAA \phi(x_\ell)  +
  \frac{h}{2} \, p_t(x_L\,\vert\,x)\, \AAA \phi(x_L) \ \\
& \simeq \sum_{\ell = 0}^L P_t(l) \,  \AAA \phi(x_\ell)
\end{align*} 
with
\begin{equation}
\label{eq.Ptpt} 
  P_t(0) \simeq  \frac{h}{2} \, p_t(0\,\vert\,x), \quad
  P_t(L) \simeq  \frac{h}{2} \, p_t(x_L\,\vert\,x), \quad
  P_t(\ell) \simeq  h \, p_t(x_\ell\,\vert\,x), \, \text{ for } 0 < \ell < L \,.
\end{equation} 
We also need to define~$P_t(\Upsilon) \simeq E_t(x)$. When designing
our approximation, we expect~$P_t(\cdot)$ to be a discrete probability
distribution on~$\{\Upsilon, x_0, \dots, x_L\}$. The differential
operator $\AAA$ is now replaced by its finite difference
approximation, denoted $A$, using an up--wind scheme, which reads for
an interior point~$x_\ell$ with $1\leq  \ell 
\leq  L - 1$: 
\begin{align*} 
  \phi'(x_\ell) 
  &\simeq 
  \begin{dcases}
    \frac{\phi(x_{\ell+1}) - \phi(x_\ell)}{h}, \quad \text{if } b(x_\ell) \geq
    0\,,\\
    \frac{\phi(x_\ell) - \phi(x_{\ell-1})}{h}, \quad \text{if } b(x_\ell) < 0\,,\\
  \end{dcases} \\[0.5em]
  \phi''(x_\ell) 
  &\simeq 
  \frac{\phi(x_{\ell+1}) - 2 \, \phi(x_\ell) + \phi(x_{\ell-1}) }{h ^2 }\,.
\end{align*} 
The resulting approximation can be written as
\begin{displaymath}
  \AAA\, \phi(x_\ell) \simeq A_{\ell,\ell-1}\, \phi(x_{\ell-1}) +
  A_{\ell,\ell}\, \phi(x_\ell) + 
  A_{\ell,\ell+1}\, \phi(x_{\ell+1}), \quad \ell = 1,\dots,(L-1) \, , 
\end{displaymath} 
with
\begin{align*} 
\forall \, \ell = 1,\dots , L-1, \quad 
A_{\ell,\ell-1} &= \frac{b^-(x_\ell)}{h} + \frac{a(x_\ell)}{2\,
  h ^2}\,, \\
A_{\ell,\ell} &= - \frac{ | b(x_\ell)| }{h} - \frac{a(x_\ell)}{
  h ^2}\,, \\
A_{\ell,\ell+1} &= \frac{b^+(x_\ell)}{h} + \frac{a(x_\ell)}{2\,
  h ^2} \,.
\end{align*} 
Appropriate boundary condition at $x_0$ and $x_L$ will be given later
on. It is enlightening to interpret this operator $A$ as the
infinitesimal generator of a pure jump Markov process on the grid
$(\Upsilon, x_0, \dots, x_L)$. Indeed,  the  extra--diagonal terms of
$A$, considered as a matrix, are non--negative and the sum on each row
is~$0$. $P_t(\ell)$ is then the probability that this process occupies
site $x_\ell$ at time~$t$. From an interior point $x_\ell$, this
process jumps to one of its neighbors with a bias directed according
to the drift. This interpretation suggests how to complete the three
lines of~$A$ not yet defined.  We set all coefficient of the first line
to~$0$, since it corresponds to the absorbing
state~$\Upsilon$. 
We introduce the notation $\mathbf{P}_t = \left (P_t(\ell)\right
)_{\ell = \Upsilon,0,\dots ,L}$ for the law of the jump process at
time~$t$ starting from~$x$.  This probability distribution solves the
Fokker--Planck equation for jump processes that reads
\begin{equation} 
\label{eq.kf.jump}
  \dot{\mathbf P}_t = A^\ast\mathbf{P}_t \ .
\end{equation} 
Observe that the first ODE of system~\eqref{eq.kf.jump} is 
\begin{displaymath}
 \dot P_t(0) = A_{0,0} \, P_t(0) +  A_{1,0} \, P_t(1)
\end{displaymath} 
where
\begin{displaymath}
  A_{1,0} = \frac{b^-(h)}{h} + \frac{a(h)}{2\,
  h ^2}  \,.
\end{displaymath} 
Using~\eqref{eq.Ptpt}, this gives an approximation 
\begin{align} 
\label{eq.p0dot.fd} 
  \frac{\partial p_t(0\,\vert\,x)}{\partial t } 
  \simeq
 \left (2\, \frac{b^-(h)}{h} + \frac{a(h)}{h ^2} \right ) 
 \, p_t(h\,\vert\,x) +  A_{0,0}\,  p_t(0\,\vert\,x)\,.
\end{align} 
A suggestion is to find~$A_{0,0}$ such that~\eqref{eq.p0dot.fd} is a
finite difference approximation for $\lim_{y\downarrow 0} \AAA ^\ast p_t(y\,\vert\,x)$: 
\begin{align*} 
  \lim_{y\downarrow 0} \AAA ^\ast p_t(y\,\vert\,x)  
  & =  
  - b'(0)\,p_t(0\,\vert\,x) 
  + \frac12 \,a''(0)\, p_t(0\,\vert\,x) 
+   \frac{\partial p_t(y\,\vert\,x)}{\partial y }{\Big |_{y = 0}}\ 
  a'(0)  \,.
\end{align*} 
This limit involves only the first derivative of~$p_t(y\,\vert\,x)$ due to the
vanishing diffusion. With
\begin{displaymath}
  \frac{\partial p_t(y\,\vert\,x)}{\partial y } {\Big |_{y = 0}} 
  \simeq
  \frac{p_t(h\,\vert\,x) - p_t(0\,\vert\,x)}{h}
\end{displaymath} 
we obtain the approximation  
\begin{align}
\label{eq.p0dot} 
  \frac{\partial p_t(0\,\vert\,x)}{\partial t}  
  &\simeq
  p_t(0\,\vert\,x) \, \left [ 
    -  b'(0) 
    + \frac12 \, a''(0) 
    - \frac{1}{h}\,a'(0) 
  \right ] 
 +  p_t(h\,\vert\,x)  \frac{1}{h}\, a'(0) \,.
\end{align} 
Using
\begin{displaymath}
  a(h) = h\,a'(0)+h^2\, \frac12 \,a''(0) 
\end{displaymath} 
in~\eqref{eq.p0dot.fd} we get
\begin{displaymath} 
  \frac{\partial p_t(0\,\vert\,x)}{\partial t}  \simeq
  2\, \frac{b^-(h)}{h} \,  p_t(h\,\vert\,x) 
  + p_t(h\,\vert\,x)  \frac{1}{h}\, a'(0)
  + \frac12\, a''(0)\,  p_t(h\,\vert\,x) 
  +  A_{0,0}\,  p_t(0\,\vert\,x) \ .
\end{displaymath} 
Also,~$p_t(h\,\vert\,x)\simeq p_t(0\,\vert\,x)$ gives
\begin{displaymath}
  \frac{\partial p_t(0\,\vert\,x)}{\partial t}  \simeq
  p_t(0\,\vert\,x) \, \left [ 2\, \frac{b^-(h)}{h}
    + \frac12\, a''(0)
    + A_{0,0}
    \right ]
  + p_t(h\,\vert\,x)  \frac{1}{h}\, a'(0) \ .
\end{displaymath} 
Now since
\begin{displaymath}
  b^-(h) = \frac{\vert b(h) \vert - b(h) }{2}
  \quad \text{and} \quad
  b(h) \simeq h\, b'(0) 
\end{displaymath} 
we finally have
\begin{displaymath} 
\frac{\partial p_t(0\,\vert\,x)}{\partial t}  \simeq
  p_t(0\,\vert\,x) \, \left [
    |b'(0)| - b'(0) + \frac12\, a''(0)
    + A_{0,0}
    \right ]
  + p_t(h\,\vert\,x)  \frac{1}{h}\, a'(0) \ .
\end{displaymath} 
In order to have an approximation of~\eqref{eq.p0dot}, we must set 
\begin{displaymath} 
  A_{0,0} = - |b'(0)|  - \frac{1}{h} \,a'(0)\,. 
\end{displaymath} 
This diagonal term of~$A$ is non--negative as expected. We see that
the state~$0$ of the jump process is not absorbing since $A_{0,0} \not
= 0$, but act as a transition state towards
extinction~$\Upsilon$. Since there is no reason to allow a jump to an
interior point, we also set 
\begin{displaymath}
\forall \, \ell = 1,\dots L\quad  A_{0,\,\ell} = 0,\quad
\text{ and } \quad A_{0,\Upsilon} = - A_{0,0} \ . 
\end{displaymath} 
Observe that, from~\eqref{eq.kf.jump}  the probability of extinction
$P_t(\Upsilon)$ satisfies the evolution equation 
\begin{displaymath}
  \dot  P_t(\Upsilon) 
  = 
  -  A_{0,0} \, P_t(0)  
  = 
  \Big( 
     \frac12 \, a'(0)+ \frac{h}{2} \,   |b'(0)|
  \Big ) \; p_t(0\,\vert\,x) \,.
\end{displaymath} 
When $h\downarrow0$, this equation is consistent with~\eqref{eq.fpe.E}
which gives the rate of extinction. Notice that $A_{0,0}$ could have
been chosen so that the above equation exactly
matches~\eqref{eq.fpe.E}, but then~\eqref{eq.p0dot.fd}
and~\eqref{eq.p0dot} would not match so closely.

The right boundary is simple: in order for the jump process to remain
on the grid, its behavior at  boundary $x_L$ has to be prescribed 
artificially. There is no canonical choice between absorbing or
reflecting boundary condition, since both corrupt the theoretical
behavior. We choose the reflecting boundary condition at~$x_L$ that
reads:
\begin{align*}
 A_{L,L-1} &= \frac{ |b(x_L)| }{h} + \frac{a(x_L)}{
  h ^2}\,,
  &
A_{L,L} &= -\frac{|b(x_L)|}{h} - \frac{a(x_L)}{h ^2} \,.
\end{align*} 
The sum on the last row is~$0$, so that there is no probability leak
at boundary~$x_L$. The boundary condition at~$0$ requires more
care. 
\subsubsection*{Time discretization}

\begin{figure}
  \centering
  \def\svgwidth{9cm}
 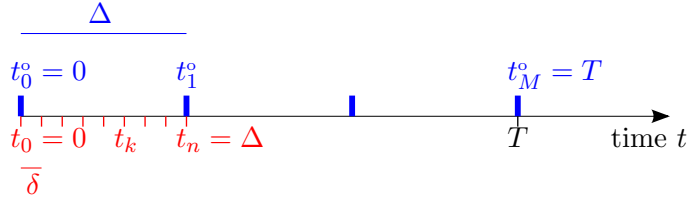
\caption{The observation instants are $\tto_{k}=k\,\Delta$ with 
$\Delta=\frac{T}{M}$; the time-discretization instants are
$t_{k}=k\,\delta$ with 
$\delta=\frac{\Delta}{n}$ for $n$ given.}
\label{fig.time}
\end{figure}

Equation~\eqref{eq.kf.jump} is discretized in time using the Euler
implicit scheme 
\begin{displaymath} 
  [I - \delta\,A]^\ast\, \mathbf{\tilde P}_{t_{k+1}} = \mathbf{\tilde P}_{t_{k}}\,,
  \quad k=0,\dots,n-1
\end{displaymath} 
where $t_k \eqdef  k\,\delta$ with $\delta=\Delta/n$, $n$ given; see
Figure \ref{fig.time}. 
The initial condition is approximated by
\begin{displaymath}
  \tilde P_{t_0}(l) = \begin {dcases}
             1\,,& \text{ if } \ell = \ell_0\,,\\
             0\,,& \text{ otherwise, }
           \end {dcases} 
\end{displaymath}  
where $x_{\ell_0}$ is the nearest neighbor in the grid of the initial
condition $x$. According to~\eqref{eq.Ptpt}, the numerical solution
$\mathbf{\tilde P}_\Delta$ 
yields a numerical approximation $\tilde p_\Delta(x_\ell\,\vert\,x)$ for the
density at a grid point, that can be linearly interpolated to obtain
an approximation $\tilde p_\Delta(y\,\vert\,x)$ for $0\leq y \leq x_L$.
The likelihood function is then approximated by
\begin{displaymath}
q_\Delta (x,y) \simeq
  \begin {dcases}
    \textstyle
    P_\Delta(\Upsilon)\,,
    &
    \text{ if } y = 0\,,
\\
   \tilde p_\Delta(y\,\vert\,x)\,,
   & \textrm{ if } y \in ]0,x_L] \,.
  \end {dcases} 
\end{displaymath} 
%


\begin{remark} 
This discretization scheme is unconditionally stable, but $h$ and
$\delta$ have to be chosen in a coherent way. Indeed, 
$-A(\ell,\ell)$ gives the expectation of the holding time of the pure
jump Markov process. We see that the order of magnitude of the holding
time is $\frac{1}{h ^2}$. The time step~$\delta$ should then be chosen small
enough to ensure that not too many jumps occur within an interval of
length~$\delta$.
\end{remark} 
%
The numerical treatment of the Fokker--Planck equation in the
degenerate case has already been considered in the numerical analysis,
see for example~\citet{cacio2011feller}. The
approach adopted in this work retains the
probabilistic meaning of the objects involved, at the cost of a
possible loss of accuracy.

\subsection{Monte Carlo approximations}
\label{sec.approximation.MC}
%
A number of estimation methods using Monte Carlo simulations have been  
proposed in the absolutely continuous case, see~\citet{hurn2007a} and
references therein for a detailed 
review, or~\citet{fearnhead2008a} which includes many application
examples. In this section, we will design 
modifications of some of them 
in order to handle the extinction feature.

The numerical approximations of Monte Carlo type presented hereafter
involve simulation of a $N$--samples with common
law~$Q_\Delta^\theta(\rmd y\,\vert\,x)$, for many different initial
conditions~$x$. In our case, we will not be able to draw random
variates from~$Q_\Delta^\theta( \rmd y\,\vert\,x)$ exactly, but only
from distributions close to it. The simplest algorithm for simulating
trajectories of~\eqref{eq.sde} is the Euler--Maruyama scheme,
restricted to non--negative values, that is
\begin{align}
\label{eq.euler.positif} 
  \bar X_{t_{k+1}}  
  &= 
  \max\Big(
  0
  , 
  \bar X_{t_{k}} + \delta\,  b(\bar X_{t_{k}})
   + \sqrt{\delta} \, \sigma (\bar X_{t_{k}})\, w_{k}
  \Big) 
  \,,\quad k=0,\dots,n-1
\end{align} 
with $\bar X_{0}  = x$ and where $w_{k}$ are i.i.d. $\NN(0,1)$. One
iteration of the approximation scheme \eqref{eq.euler.positif} amounts
to draw from the  transition kernel $K_\delta(\rmd z\,\vert\,x)$ instead
of~$Q_\delta(\rmd z\,\vert\,x)$ where $K_\delta(\rmd z\,\vert\,x)$ is
defined  by:
\begin{equation}
\label{eq.euler.positif.kernel} 
  K_\delta(\rmd z\,\vert\,x) 
  \eqdef 
  \begin{cases}
    e_{\delta}(x) \, \delta_0(\rmd z)   + g_\delta(z\,\vert\,x)\,\rmd z\,,
    &\textrm{ if }x>0\,,
    \\
    \delta_0(\rmd z)\,,
    &\textrm{ if }x=0\,,
  \end{cases}
\end{equation} 
with
\begin{align*}
  e_{\delta}(x)
  &\eqdef
  \textstyle
  1 -  \int_0 ^\infty g_\delta(z\,\vert\,x)\, \rmd z\,,
\\
  g_\delta(z\,\vert\,x) 
  &\eqdef 
  \textstyle
     \frac{1}{ \sqrt{2\, \pi \, \delta \,\sigma (x)}} \, 
     \exp \left \{-  \frac{(z - x - \delta\, b(x))^2}{2\, \delta \,
       \sigma(x)} \right \}\, \mathbf{1}_{\R_+}(z)  \,.
\end{align*} 
It is known that this approximation is valid if $\delta$ is
sufficiently small, see~\citet[chapter~4]{risken1996a}. 
The numerical scheme~\eqref{eq.euler.positif} produces samples of
common law~$\QE_\Delta(\rmd y\,\vert\,x)$
\begin{equation}
\label{eq.Qtilde} 
     \QE_\Delta(\rmd y\,\vert\,x)  
     =  
     \int_{\{y_1\geq 0\}} \cdots  \int_{\{y_{n-1} \geq  0\}} 
         K_\delta (\rmd y_1\,\vert\,x)\,   K_\delta(\rmd
         y_2\,\vert\,y_1)\, \cdots       
          K_\delta(\rmd y\,\vert\,y_{n-1}) \ ,
\end{equation} 
which is an approximation of the true transition
kernel~$Q_\Delta(\rmd y\,\vert\,x)$. Notice that, using the semi--group and
the Markov properties, we also have a similar
decomposition
\begin{displaymath}
   Q_\Delta(\rmd y\,\vert\,x) 
   = 
   \int_{\{y_1\geq 0\}} \cdots\int_{\{y_{n-1} \geq 0\}}  
       Q_\delta (\rmd y_1\,\vert\,x)\, Q_\delta(\rmd y_2\,\vert\,y_1)\, \cdots
            Q_\delta(\rmd y\,\vert\,y_{n-1}) \ .
\end{displaymath}
\begin{remark} 
The recent works~\citet{beskos2005a} about the
\emph{Exact Algorithm}~(EA) seems promising for drawing exactly
from~$Q_\Delta(\rmd y\,\vert\,x)$. To our knowledge, this algorithm cannot be
applied directly to our specific case, due to the almost sure
extinction. There are also other alternatives
to~\eqref{eq.euler.positif}, such as the Milstein scheme, see
e.g.~\citet{kloeden2003a}, or Euler--Maruyama scheme for killed
diffusion, see~\citet{gobet2001a}. Modifications of these
algorithms might be necessary to cope with the extinction
problem.
\end{remark} 
\begin{remark} 
If $\Delta$ is itself 
small enough, there would be no need to simulate the solution
of~\eqref{eq.sde} at intermediate time between $0$ and~$\Delta$, since
$K_\Delta$ has an explicit density with respect to~$m$. However, for
most applications, the observations are not available at a so high
sampling rate. 
\end{remark} 

\subsubsection{Non--parametric estimation}

A simple usage of the approximation scheme~\eqref{eq.euler.positif}
is to produce a $N$--sample $\bar X_{\Delta}^{(1)},\dots \bar
X_{\Delta}^{(N)}$ of $\QE_\Delta(\rmd y\,\vert\,x)$. These faked
observations are then fed into a nonparametric 
estimate of the density
$p_\Delta(y\,\vert\,x)$ at~$y$, denoted by $\hat
p_\Delta(y\,\vert\,x)$. Again, the case of extinction should be cared
for by first discarding the values $\bar X_{\Delta}^{(i)} = 0$, if
there are any, from the estimation. The resulting approximation reads
\begin{displaymath}
  q_\Delta (y\,\vert\,x) 
  \simeq
  \begin {cases}
    \frac{N - N_s}{N} \,,   &\text{ if } y = 0\,,\\
    \frac{N_s}{N} \,\hat p_\Delta (y\,\vert\,x)\,, & \text{ otherwise } 
  \end {cases} 
\end{displaymath} 
where $N_s$ is the random number of trajectories still alive at time
$\Delta$, that is
\[
 N_s \eqdef \#\big\{i=1,\dots,N\,;\, \bar X_{\Delta}^{(i)} \ne 0\big\}\,.
\]
This approach, without extinction,  is presented
in~\citet{hurn2003a}. Its efficiency relies on that of the  
nonparametric estimation method used and is therefore subject to the
classical problems of choice of bandwidth and leakage of mass in
inaccessible region ($\R_-$ in our case).

\subsubsection{Pedersen method} 
\label{sec.pedersen}
\begin{figure}
  \centering
  \def\svgwidth{10cm}
  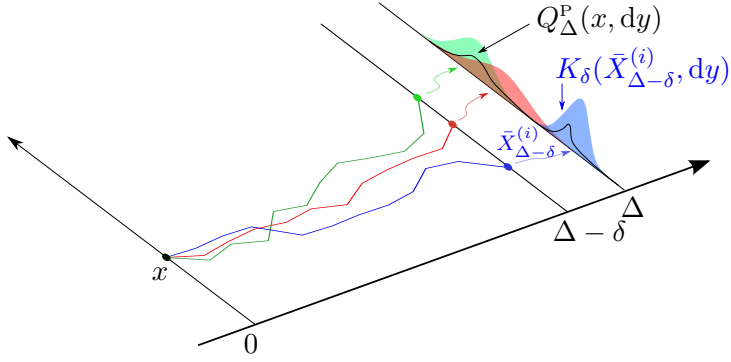
\caption{The Perdersen approximation $\QP_\Delta$ of the kernel 
$Q_\Delta$ is obtained by sampling $N$ independent trajectories $\bar
  X_{t_{0:n-1}}^{(i)}$ from \eqref{eq.euler.positif}; then the
  approximation is given by \eqref{eq.Q.pedersen}. See text for the
  precise treatment of extinct trajectories.} 
\label{fig.pedersen} 
\end{figure}

It is possible to avoid the non--parametric estimation stage. Indeed from
the Markov property:
\begin{equation}
\label{eq.markov} 
  Q_\Delta(\rmd y\,\vert\,x)  
  = 
  (Q_{\Delta - \delta}Q_\delta)(\rmd y\,\vert\,x) 
  = 
  \int_{\{y_{n-1} \geq 0\}}   
    Q_{\Delta - \delta} (\rmd y_{n-1}\,\vert\,x)\, Q_\delta(\rmd y\,\vert\,x) \,,
\end{equation} 
first approximate $ Q_{\Delta - \delta} (\rmd y_{n-1}\,\vert\,x)$ by 
  $\QE_{\Delta - \delta} (\rmd y_{n-1}\,\vert\,x)$, hence
\begin{displaymath}
   Q_{\Delta - \delta} (\rmd y_{n-1}\,\vert\,x) 
   \simeq 
   \frac{1}{N}\,\sum_{i=1}^N \delta_{\bar X_{t_{n-1}}^{(i)}}(\rmd y_{n-1}) 
\end{displaymath} 
where $\bar X_{t_{n-1}}^{(i)} \simiid  \QE_{\Delta- \delta}(\rmd y_{n-1}\,\vert\,x)$, $i=1,\dots,N$;
then approximating $Q_\delta(\rmd y\,\vert\,y_{n-1})$ by $K_\delta(\rmd
y\,\vert\,y_{n-1})$, leads to the following approximation of the
kernel $Q_\Delta(\rmd y\,\vert\,x)$: 
\begin{align}
\label{eq.Q.pedersen}
 \QP_\Delta(\rmd y\,\vert\,x)  
 &\eqdef
 \frac{1}{N}\,\sum_{i=1}^N K_\delta(\rmd y\,\vert\,\bar X_{t_{n-1}}^{(i)})  
\end{align} 
see Figure \ref{fig.pedersen}.
Let us re--numbered the sampled trajectories so that the surviving ones
correspond to $i=1,\dots,N_{s}$, according to
\eqref{eq.euler.positif.kernel} we get: 
\begin{align*}
 \QP_\Delta(\rmd y\,\vert\,x)  
 &= 
 \textstyle
 \frac{1}{N}\,\sum_{i=N_{s}+1}^N K_\delta(\rmd y\,\vert\,0)  
 +
 \frac{1}{N}\,\sum_{i=1}^{N_{s}} K_\delta(\rmd y\,\vert\,\bar X_{t_{n-1}}^{(i)})  
\\
 &= 
 \textstyle
 \frac{N - N_s}{N} \, \delta_0(\rmd y)   
      + \frac{1}{N}\, \sum_{i=1}^{N_s} 
      \big[  
        e_{\delta}(\bar X_{t_{n-1}}^{(i)})\,\delta_{0}(\rmd y)  
        +
        g_{\delta}(y\,\vert\, \bar X_{t_{n-1}}^{(i)})\,\rmd y
      \big]
\\
 &= 
 \textstyle
 \big[
    \frac{N - N_s}{N}
    +
    \frac{1}{N}\, \sum_{i=1}^{N_s}e_{\delta}(\bar X_{t_{n-1}}^{(i)})
 \big] 
 \, \delta_0(\rmd y)   
  + 
  \frac{1}{N}\, \sum_{i=1}^{N_s} 
        g_{\delta}(y\,\vert\, \bar X_{t_{n-1}}^{(i)})\,\rmd y
\end{align*} 
so that $\QP_\Delta(\rmd y\,\vert\,x)$ admits the following density with
respect to the measure $m(\rmd y)$: 
\begin{displaymath}
  \qP_\Delta (y\,\vert\,x) 
  \eqdef
  \begin {cases}
    \displaystyle
    \frac{N - N_s}{N}
    +
    \frac{1}{N}\, \sum_{i=1}^{N_s}e_{\delta}(\bar X_{t_{n-1}}^{(i)})\,,
    & \text{ if } y = 0\,,
    \\
    \displaystyle
    \frac{1}{N}\, \sum_{i=1}^{N_s} 
        g_{\delta}(y\,\vert\, \bar X_{t_{n-1}}^{(i)})\,,
   &\text{ otherwise. }  
  \end {cases} 
\end{displaymath} 
This approach is presented in~\citet{pedersen1995new} for
diffusion on $\R^d$ having an absolutely continuous density. Even with
our  adaptation allowing the extinction, it is
very easy to implement and does not involve heavy computations. It
suffers however from a well known  problem of Monte Carlo methods: if
not so many trajectories terminate around the observation~$y$ at which
the density is to be evaluated, then only a few terms will
significantly contribute to the approximation of~$p_\Delta(\rmd
y\,\vert\,x)$. In this case the approximation will be of poor
quality. Beside, a number of trajectories will have been generated
uselessly. This problem is usually well tackled by importance sampling
procedures.
%

\subsubsection{Importance sampling with a Brownian bridge}
\label{sec.bb}

\begin{figure}
  \centering
  \def\svgwidth{10cm}
  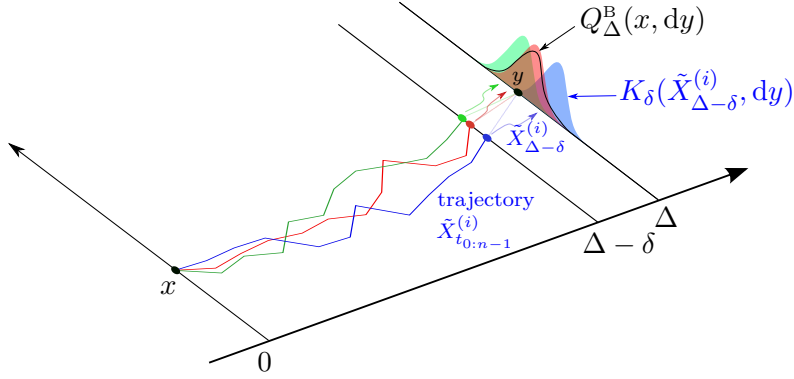
\caption{The approximation $\QB_{\Delta}$ by importance sampling with
  a Brownian bridge of the kernel $Q_{\Delta}$ is obtained by sampling
  $N$ independent trajectories $\tilde X_{t_{0:n-1}}^{(i)}$ of the
  approximation of the Brownian bridge \eqref{eq.euler.bb}; then the
  approximation is given by \eqref{eq.Q.bridge}. See text for the
  precise treatment of extinct trajectories.} 
\label{fig.bridge}
\end{figure}

When using importance sampling, the trajectories from $x$ to $y_{n-1}$
are sampled according to a different process and weighted to correct
the change of law. Such a procedure is described among others
in~\citet{durham2002numerical}. Once again, we have  
to modify the method to account for the possible extinction. 
Our choice is to generate the trajectories according to 
\begin{equation}
\label{eq.euler.bb} 
  \tilde X_{t_{k+1}}  
  = 
  \max\Big(
    0
    , 
    \tilde X_{t_{k}} + \delta\,  \frac{y - \tilde X_{t_{k}}}{\Delta - {t_{k}}} 
   + \sqrt{\delta} \, \sigma (\tilde X_{t_{k}})\, w_{k}
  \Big) 
  \,,\quad k=0,\dots,n-1
\end{equation} 
whith $X_0 = x$ and  $w_{k}\simiid \NN(0,1)$
which is nothing but a (modified for extinction) Euler--Maruyama scheme for the
SDE
\begin{equation} 
\label{eq.sde.bb}
\rmd X_t = \frac{y -  X_t}{\Delta - t} \, \rmd t + \sigma (X_t) \,
\rmd W_t, \quad 0\leq t < \Delta \,.
\end{equation} 
The drift term is designed  so as to force the trajectories towards the
given final value~$y$. Note that solution of~\eqref{eq.sde.bb} would
be a ``true'' Brownian bridge if $\sigma$ were constant. The transition
kernel associated with~\eqref{eq.euler.bb} depends on~$t<\Delta$ and
reads:
\begin{equation}
\label{eq.euler.bb.kernel} 
  \tilde K_{t,\delta}(\rmd z\,\vert\, x) 
  =
  \begin{cases} 
    \tilde e_{\delta}(x)\, \delta_0(\rmd z) + \tilde g_{t,\delta}(z\,\vert\, x)\,\rmd z  
    \,, &\textrm{if }x>0\,,
    \\
    \delta_0(\rmd z)   \,,
    &\textrm{if }x=0\,,
  \end{cases} 
\end{equation} 
where 
\begin{align*}
  \tilde e_{\delta}(x)
  &\eqdef
  \textstyle
  \left (1 - \int_0 ^\infty \tilde g_{t,\delta}(z\,\vert\, x)\, \rmd z \right )\,,
\\
  \tilde g_{t,\delta}(z\,\vert\, x) 
  &\eqdef
  \textstyle
     \frac{1}{ \sqrt{2\, \pi \, \delta \,\sigma (x)}} \, 
     \exp \Big \{
       -  \frac{(z - x - \delta\, \frac{y - x}{\Delta - t} )^2}
               {2\, \delta \, \sigma(x)} 
     \Big\}\, \mathbf{1}_{\R_+}(z)  \,.
\end{align*} 
This transition kernel is absolutely continuous with
respect to $K_\delta(\rmd z\,\vert\, x)$:
\[
\tilde  K_\delta(\rmd z\,\vert\, x) 
  = \psi_{t,\delta}(z\,\vert\,x)\, K_{t,\delta}(\rmd z\,\vert\, x)
\]
with 
\begin{align*} 
  \psi_{t,\delta}(z\,\vert\,x) & =
  \begin {dcases}
    \frac{1 - \int_0 ^\infty g_\delta(z\,\vert\,x)\, \rmd z }
         {1 - \int_0 ^\infty \tilde g_{t,\delta}(z\,\vert\,x)\, \rmd z }\,,
    & \text{ if } z = 0\,,
    \\[0.3em]
   \displaystyle \frac{\int_0 ^\infty  g_\delta(z\,\vert\,x)\, \rmd z }
                 {\int_0 ^\infty \tilde g_{t,\delta}(z\,\vert\,x)\, \rmd z } \,,
   &\text{ otherwise, }  
   \end {dcases} 
\end{align*} 
for $x>0$ and $\psi_{t,\delta}(z\,\vert\,0) = \mathbf{1}_{\R_+}(z)$.
We now have another expression for $\QE_{\Delta - \delta}(\rmd
y_{n-1}\,\vert\,x)$ as 
\begin{displaymath}
 \QE_{\Delta - \delta}(\rmd y_{n-1}\,\vert\,x) = \int_0^\infty \dots
 \int_0^\infty 
\Psi(y_{1},\dots,y_{n-1}\,\vert\,x)\,  \tilde K_{0,\delta} (\rmd y_1\,\vert\,x) \dots 
\tilde K_{\Delta - 2\,  \delta,\delta}(\rmd y_{n-1}\,\vert\,y_{n-2}) 
\end{displaymath} 
where
\begin{displaymath}
 \Psi(y_{1},\dots,y_{n-1}\,\vert\,x) = \psi_{\delta,\delta}(y_{1}\,\vert\,x) \dots \psi_{\Delta -
   2\, \delta,\delta}(y_{n-1}\,\vert\,y_{n-2}) \,.
\end{displaymath} 
Hereafter, we will denote by
\[
  \tilde X_{t_{1:n-1}}
  =
  (\tilde X_{t_{1}},\dots,\tilde X_{t_{n-1}})\in \R^n
\]
  a trajectory
generated by~\eqref{eq.euler.bb} up to time~$\Delta - \delta$, with
initial value $\tilde X_{0}=x$. 

\medskip

With this setting, and like in the previous section, in
\eqref{eq.markov}  first approximate $Q_{\Delta - \delta}$ by a
weighted sample
\begin{displaymath}
  Q_{\Delta - \delta} (\rmd y_{n-1}\,\vert\,x) 
  \simeq 
  \frac{1}{N}\, 
  \sum_{i=1}^N \Psi(\tilde X_{t_{1:n-1}}^{(i)}\,\vert\,x)\, 
     \delta_{\tilde X_{t_{n-1}}^{(i)}}(\rmd y_{n-1}) 
\end{displaymath} 
then, as before, approximate $Q_\delta(\rmd y\,\vert\,y_{n-1})$ by
$K_\delta(\rmd y\,\vert\,y_{n-1})$; finally the kernel $Q_\Delta(\rmd
y\,\vert\,x)$ defined by \eqref{eq.markov} is approximated by:  
\begin{align}
\label{eq.Q.bridge}
   \QB_\Delta(\rmd y\,\vert\,x)  
   &\eqdef 
   \frac{1}{N}\, 
   \sum_{i=1}^N \Psi(\tilde X_{t_{1:n-1}}^{(i)}\,\vert\,x)\, 
       K_\delta(\rmd  y\,\vert\,\tilde X_{t_{n-1}}^{(i)})  \,.
\end{align} 
As before, let us re--numbered the sampled trajectories so that the
surviving ones correspond to $i=1,\dots,N_{s}$:  
\begin{align*}
   \QB_\Delta(\rmd y\,\vert\,x)  
   &=  
   \textstyle
   \frac{1}{N}\, 
   \sum_{i=N_{s}+1}^N \Psi(\tilde X_{t_{1:n-1}}^{(i)}\,\vert\,x)\, 
       K_\delta(\rmd  y\,\vert\,0)  
   +
   \frac{1}{N}\, 
   \sum_{i=1}^{N_{s}} \Psi(\tilde X_{t_{1:n-1}}^{(i)}\,\vert\,x)\, 
       K_\delta(\rmd  y\,\vert\,\tilde X_{t_{n-1}}^{(i)})  
\\
   &=  
   \textstyle
   \frac{1}{N}\, 
   \sum_{i=N_{s}+1}^N \Psi(\tilde X_{t_{1:n-1}}^{(i)}\,\vert\,x)\, 
       \delta_{0}(\rmd  y)  
   \\
   &\textstyle\qquad\qquad
   +
   \frac{1}{N}\, 
   \sum_{i=1}^{N_{s}} \Psi(\tilde X_{t_{1:n-1}}^{(i)}\,\vert\,x)\, 
       \big[
       e_{\delta}(\tilde X_{t_{n-1}}^{(i)})\,\delta_{0}(\rmd  y)
       +
       g_{\delta}(y\,\vert\,\tilde X_{t_{n-1}}^{(i)})\,\rmd  y
       \big] 
\\
   &=  
   \textstyle
   \big[
   \frac{1}{N}\, 
   \sum_{i=N_{s}+1}^N \Psi(\tilde X_{t_{1:n-1}}^{(i)}\,\vert\,x)
   +
   \frac{1}{N}\, \sum_{i=1}^{N_{s}}
       \Psi(\tilde X_{t_{1:n-1}}^{(i)}\,\vert\,x)\,e_{\delta}(\tilde
       X_{t_{n-1}}^{(i)}) 
   \big]\, 
       \delta_{0}(\rmd  y)  
   \\
   &\textstyle\qquad\qquad
   +
   \frac{1}{N}\, 
   \sum_{i=1}^{N_{s}} \Psi(\tilde X_{t_{1:n-1}}^{(i)}\,\vert\,x)\, 
       g_{\delta}(y\,\vert\,\tilde X_{t_{n-1}}^{(i)})\,\rmd  y
\end{align*} 
hence $\QB_\Delta(\rmd y\,\vert\,x)$ admits the following density with
respect to  $m(\rmd y) $: 
\begin{displaymath}
  \qB_\Delta (y\,\vert\,x) 
  \eqdef
  \begin {dcases}
    \frac{1}{N}\, 
   \sum_{i=N_{s}+1}^N \Psi(\tilde X_{t_{1:n-1}}^{(i)}\,\vert\,x)
   +
   \frac{1}{N}\, \sum_{i=1}^{N_{s}}
       \Psi(\tilde X_{t_{1:n-1}}^{(i)}\,\vert\,x)\,e_{\delta}(\tilde X_{t_{n-1}}^{(i)})
       \,,
   &     
  \text{ if } y = 0\,,
  \\
    \frac{1}{N}\, 
   \sum_{i=1}^{N_{s}} \Psi(\tilde X_{t_{1:n-1}}^{(i)}\,\vert\,x)\, 
       g_{\delta}(y\,\vert\,\tilde X_{t_{n-1}}^{(i)})\,,
   &\text{ otherwise. }  
  \end {dcases} 
\end{displaymath} 
\begin{remark} 
It is possible to use other importance sampler
than~\eqref{eq.euler.bb}.~\citet{durham2002numerical} use a
\emph{modified Brownian 
  bridge} which seems to reduce further the variance of the resulting
approximation. The generating algorithm reads 
\begin{displaymath} 
 \tilde X_{t_{k+1}} 
 = 
 \tilde X_{t_{k}} 
 + \delta\,  \frac{y - \tilde X_{t_{k}}}{\Delta - t} 
 + \sqrt{\delta} \, \sigma (\tilde X_{t_{k}})\, 
 (1 - \frac{\delta}{\Delta - t})\,
 w_{k}  \,.
\end{displaymath} 
The drift term is as in~\eqref{eq.euler.bb} but the variance is
progressively damped up to final time $t = \Delta-2\, h$ at which it
equals $\frac{1}{2}\,\sigma (Y_t)$.
\end{remark}

\section{Numerical experiments}
\label{sec.numerical.experiments}
We now evaluate the numerical performance of the
approximation methods described above, on the scenario defined by the
model parameter given in Table~\ref{tab.par}. These values are chosen
so  as to observe a short transient phase where the population grows
quickly, starting from a small initial value, and a stationary phase
with a high noise intensity. The simulation parameters of the
approximation methods are also given by Table~\ref{tab.par}. In
particular, the same time step~$\delta$ is used for all methods.  
\begin{table}
  \begin{center} 
    \begin{tabular}{c|c|c|c|c|c|c} 
      \multicolumn{7}{c}{Model}\\ \hline
      $\lambda$ &$\mu$ & $\alpha$ & $\rho$ & $x_0$ & $T$ & $N$ \\ \hline \hline
      $20$ & $18$ & $1$ & $10^{-1}$ & 0.25 & 10 & 200
    \end{tabular} 
\qquad
    \begin{tabular}{c|c|c} 
      \multicolumn{3}{c}{Simulation}\\ \hline
      $h$      & $\delta$ & $N$  \\ \hline \hline
      $10^{-3}$ & $10^{-3}$ & $200$ 
    \end{tabular} 
  \end{center}  
  \caption{Parameters}
  \label{tab.par}
\end{table} 
%

\subsection{Dynamics}

\begin{figure}
\begin{center}
\includegraphics[width=0.48\linewidth]{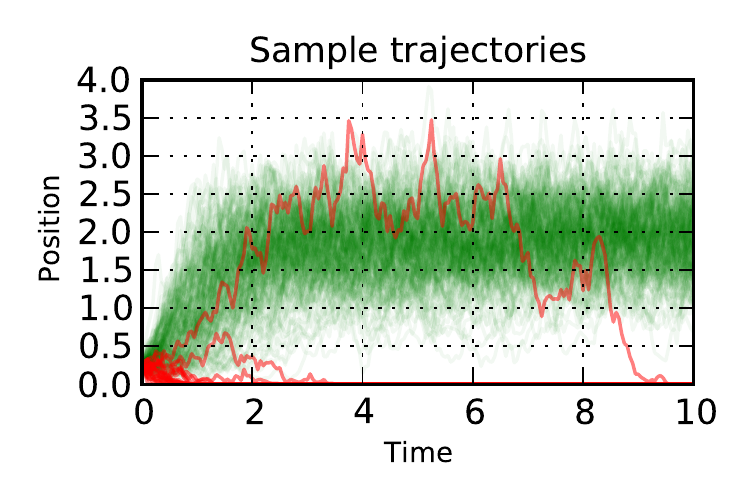}
\includegraphics[width=0.48\linewidth]{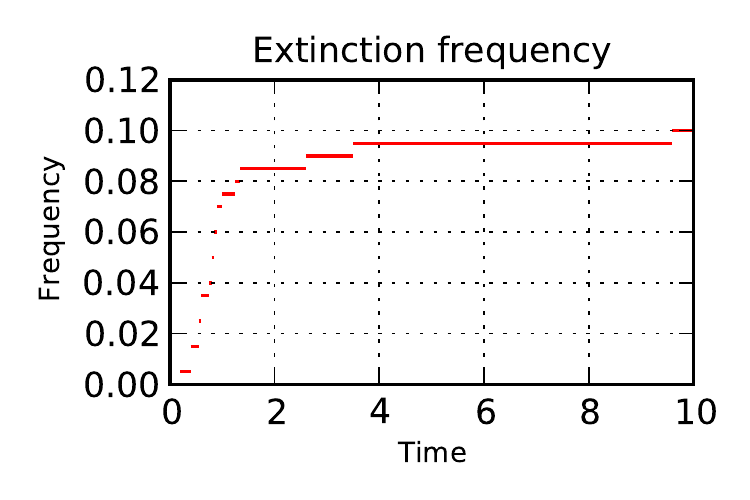}
\end{center}
\caption{\it On the left: simulation of $n=200$ trajectories according to
  dynamics~\eqref{eq.sde}. Trajectories in red become extinct before
  final time. One of them completes the transient phase. On the right
  plot, the corresponding extinction frequency.}
\label{fig.mc1}
\end{figure}

Figure~\ref{fig.mc1} shows~200 trajectories of
model~\eqref{eq.sde}, simulated with the Euler--Maruyama scheme. Within the time
interval of the simulation, we observe the three possible behaviors:
\begin{itemize}
\item Early extinction during the transient phase. This stems from the
  fact that the initial level of the population is small and the noise
  intensity is high. 
\item Completed transient phase followed by noisy fluctuations around
  a natural carrying capacity. Even if it is known that all
  trajectories will eventually reach~0, all of them but one survive on
  the short term. 
\item Late extinction after reaching the stationary phase. Only one
  trajectory is concerned.
\end{itemize} 

\begin{figure}
\begin{center}
\includegraphics[width=0.48\linewidth]{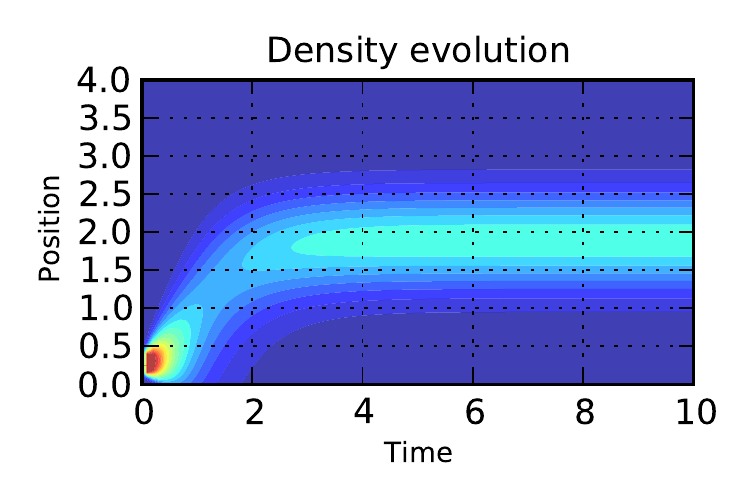}
\includegraphics[width=0.48\linewidth]{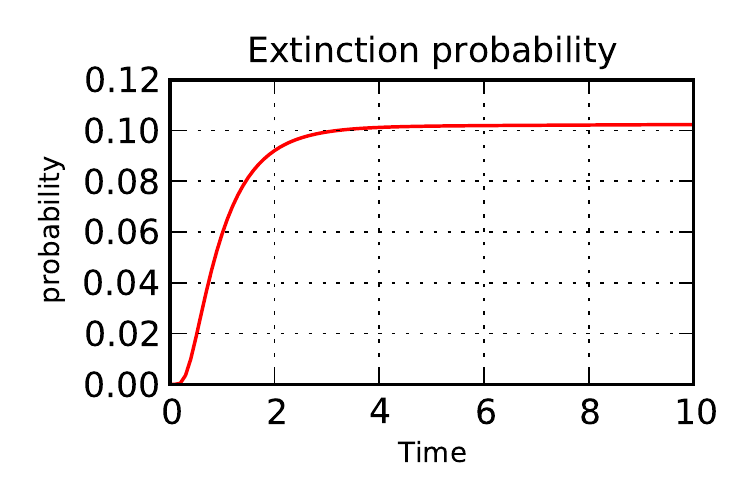}
\end{center}
\caption{\it Finite difference approximation of Fokker--Planck
  equation~\eqref{eq.fpe}. Diffusion induces a quick loss of mass in
  the transient phase. The loss of mass is then much slower.} 
\label{fig.fpe1}
\end{figure}

Due to early extinction, the estimated extinction probability
grows quickly in the first part, and seems to approach an asymptotic
value. Running the simulation on a much larger period of time
will result in an extinction frequency reaching~$1$ slowly. This
Figure is to be compared with Figure~\ref{fig.fpe1} showing a contour
plot of the numerical solution of Fokker--Planck
equation~\eqref{eq.fpe} obtained by implementing the finite difference
scheme of section~\ref{sec.approximation.fpe}. As expected, the first
two moments quickly move to stationary values, after a transient phase
where a large amount of mass is lost. On the long term, the mass keeps
on decreasing slowly.

\subsection{Kernel approximations}

\begin{figure}
\begin{center}
\includegraphics[width=\linewidth]{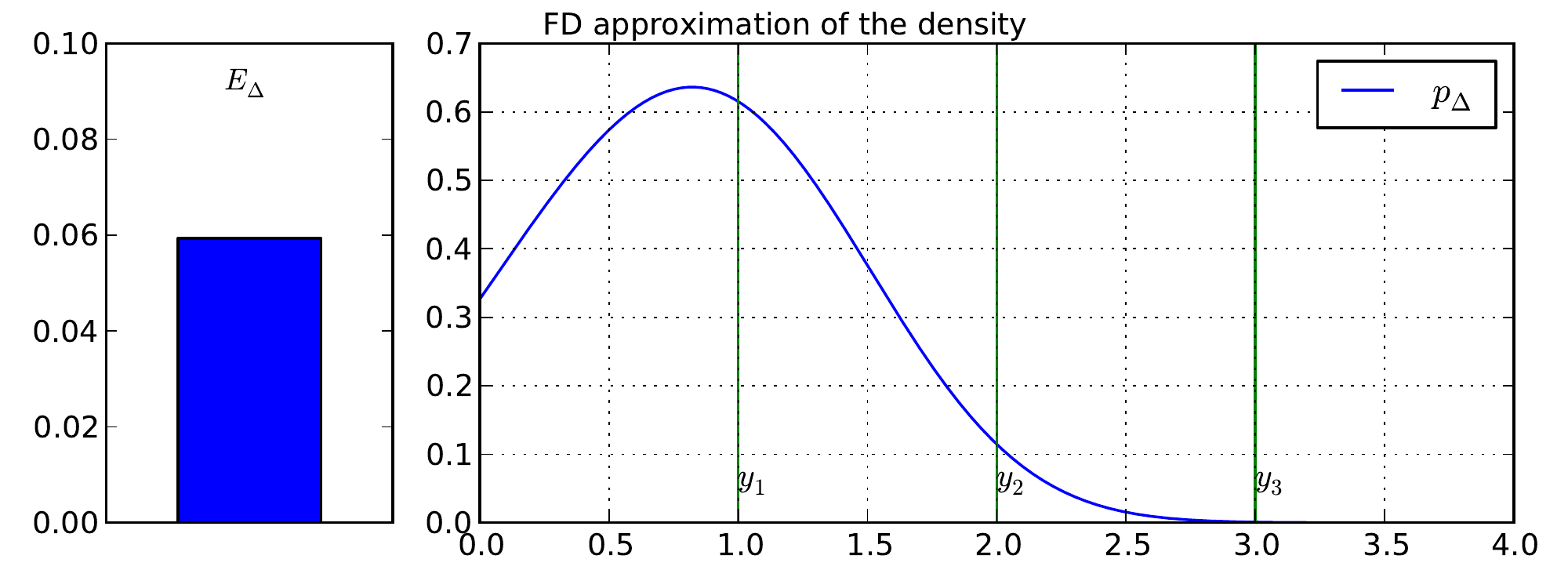} \\
\includegraphics[width=\linewidth]{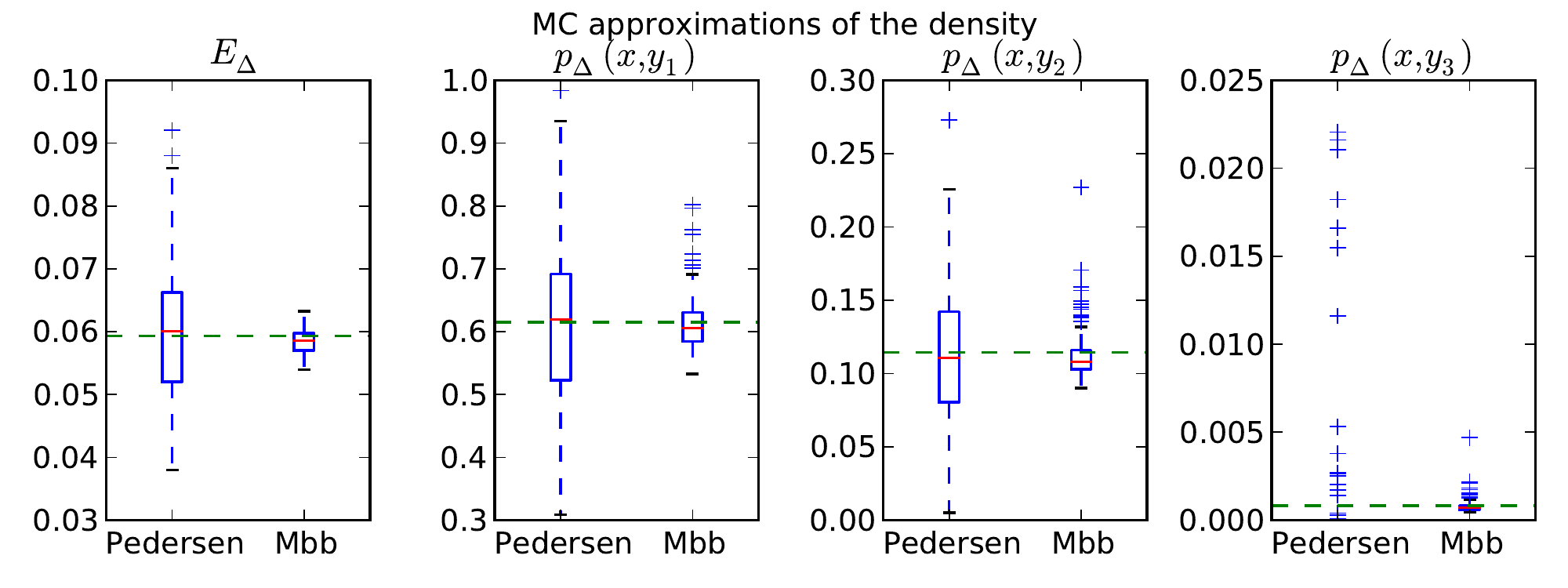}
\end{center} 
\caption{\it Comparison of approximation methods for the approximation
  of the density which solves equation~\eqref{eq.fpe}. The top row
  shows the finite difference approximation of the density~$q_\Delta$. The
  left plot is the discrete component and the right plot shows the
  continuous one. The four plots on the bottom rows give the Monte Carlo
  approximations of the density at four different locations. The green
  horizontal dashed line is the value given by the finite
  difference approximation. The model parameters are
  given in Table~\ref{tab.par}, except for $\Delta = 1$. The
  box--plots are based on~$200$ realizations. Both Monte Carlo
  methods use $N = 500$ 
  particles.}  
\label{fig.kernel}
\end{figure}

We now investigate the ability of the methods to approximate the
density~$q_\Delta$, since the final parameter estimation heavily
relies on the quality of this approximation. For this section, we take
$T = 1$ in order to observe correctly both continuous and
discrete parts of~$q_\Delta$. The other parameters are unchanged. The
results obtained with the finite difference method are shown on the
top row of Figure~\ref{fig.kernel}. There is no way to evaluate its
performance since the exact solution is not available. However, the
picture presented here does not change significantly when we use
smaller discretization steps.

To compare with the Monte Carlo methods, we consider the extinction
probability~$E_\Delta$ together with the value of the continuous
part~$p_\Delta(x,y)$ at three locations~$y_1, \, y_2$ and $y_3$. 
The bottom line of Figure~\ref{fig.kernel} shows the results of
200~independent realizations of the two Monte Carlo methods.
The green dashed line gives the values obtained by finite difference
approximation. As expected, the variance observed for the modified
Brownian bridge is much smaller than for the Pedersen method. However
the computational cost needed to achieve this performance is
high. It is worth remembering also that the $N$ trajectories simulated
for the Pedersen method do not depend of the value at which the 
density is evaluated. A single run of a~$N$ path allows the estimation
of the density at any point. This is not the case for the modified
Brownian bridge method for which a new run is needed for each
terminal value. Finally, we note that, even for the modified Brownian
bridge method, we observe a large number of outliers which are
likely to corrupt the final evaluation of the likelihood.

\subsection{Likelihood approximations}
\begin{figure}
\begin{center}
\includegraphics[width=\linewidth]{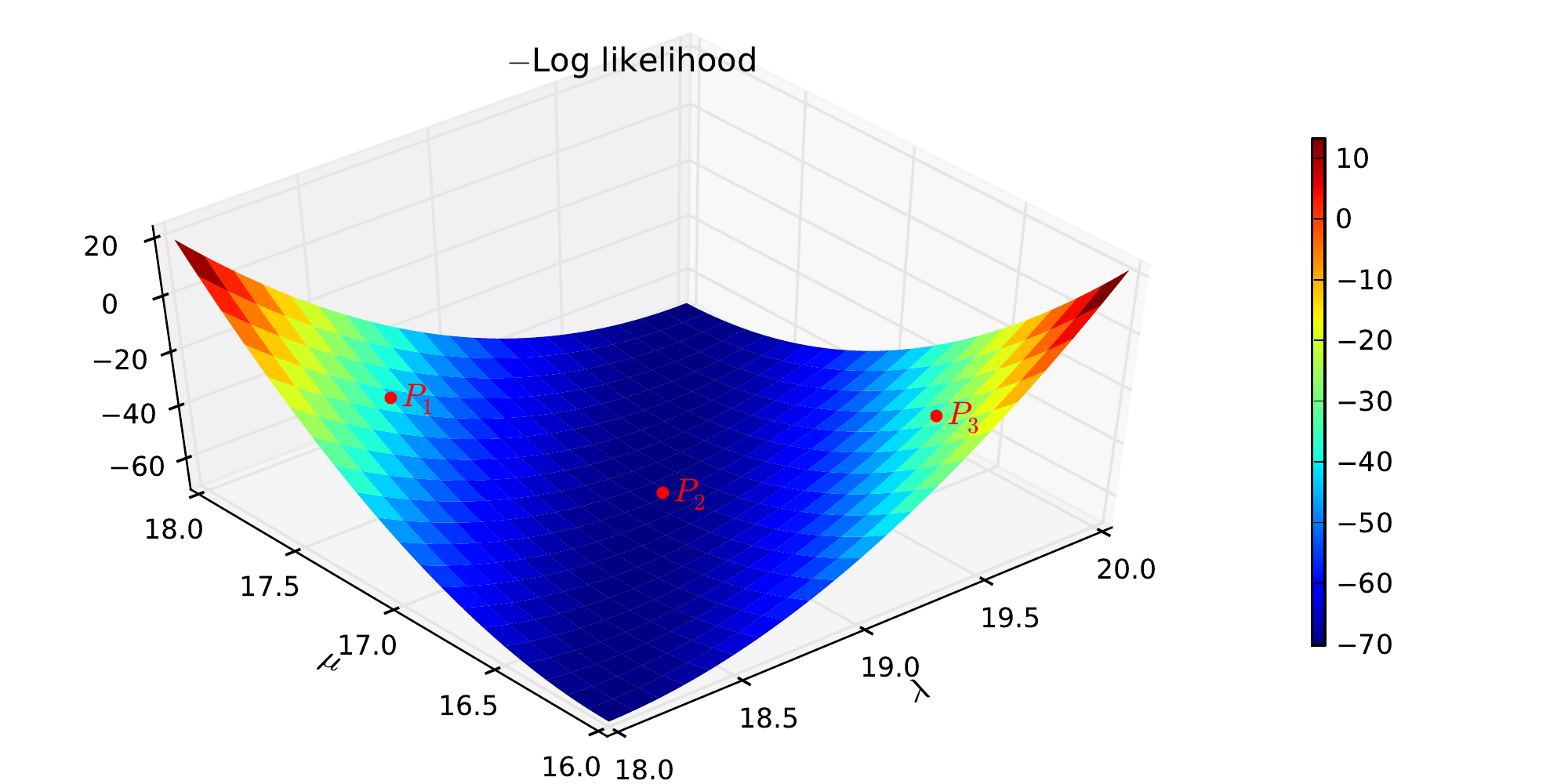} 
\end{center} 
\caption{\it Finite difference approximation of $\ell(\theta)$. The
 model parameters are given in Table~\ref{tab.par}. Initial condition
  is a Dirac mass at~$x_0= 0.25$. The  Monte Carlo approximations
  will be evaluated at points $P_1,\, P_2$   and $P_3$. }  
\label{fig.likelihoodfd}
\end{figure}

\begin{figure}
\begin{center}
\includegraphics[width=\linewidth]{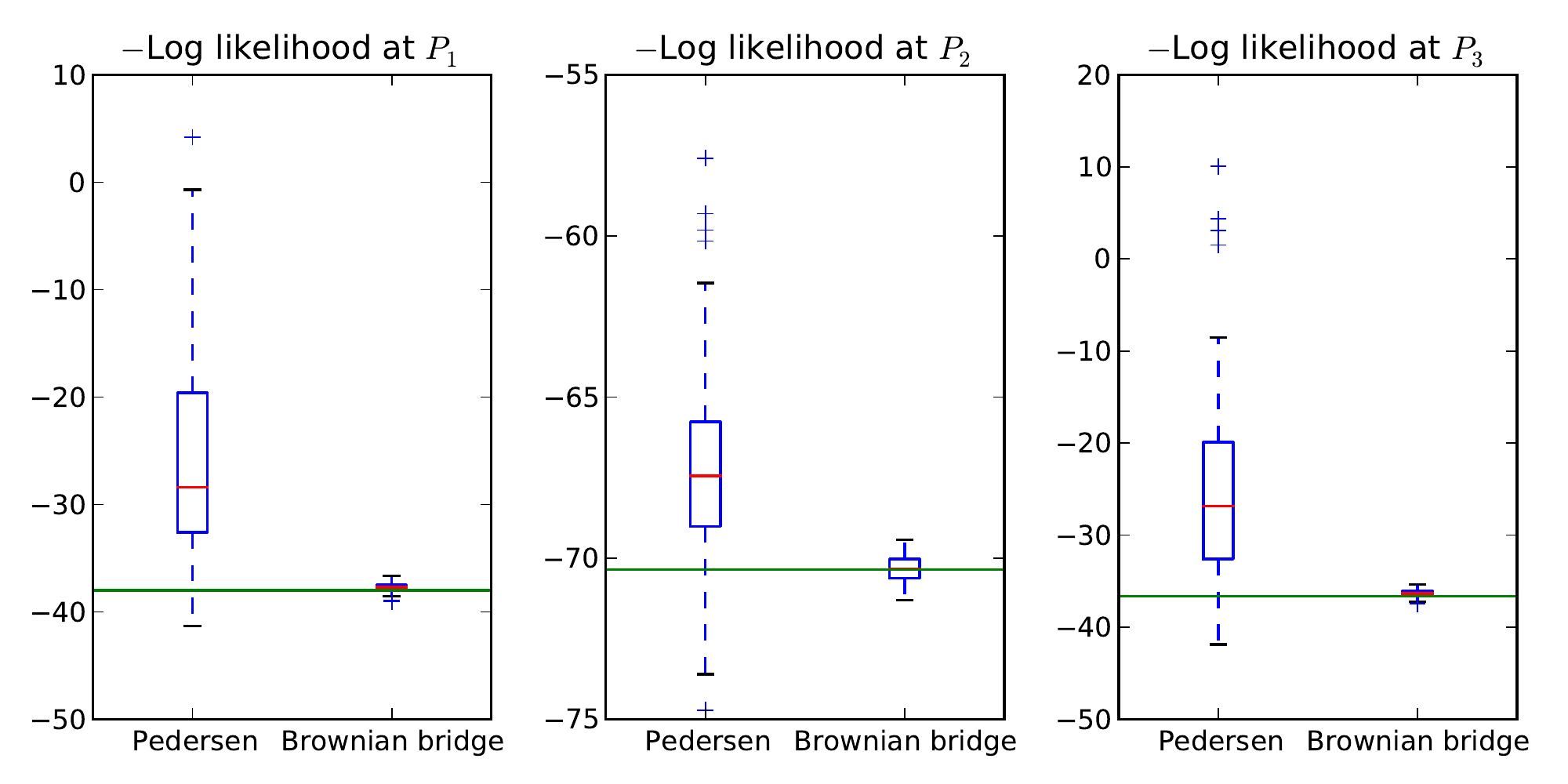}
\end{center} 
\caption{\it Comparison of Monte Carlo approximations of
  $\ell(\theta)$ at point $P_1$, $P_2$ and $P_3$ of
  Figure~\ref{fig.likelihoodfd}. The box--plots are based   on~$200$
  realizations. Both Monte Carlo methods use $N = 500$ particles.}   
\label{fig.likelihoodmc}
\end{figure}

%
The numerical results presented in this section are based on a typical
trajectory of Figure~\ref{fig.mc1}. We obtain the set of data by
sampling the trajectory at instants $k\, \Delta$, for $k =
0,\dots,200$ and $\Delta= 0.05$. We take $\theta = (\lambda,\mu) \in
\R^2_+$ as our parameter to be estimated and the other model
parameters $\alpha$, $\kappa$ as fixed known values. The initial
condition $x_0$ is also deterministic and known. To improve numerical
stability, we will note~$\ell(\theta) = - \log(\LL(\theta))$ and study
the minimum of $\ell(\theta)$. 
Figure~\ref{fig.likelihoodfd} shows the graph of the
finite difference approximation of the function to minimize
$\ell(\theta)$. We clearly distinguish two orthogonal directions of
variations defined respectively by the equations
\begin{displaymath}
\lambda + \mu = b, \quad \text{ and } \quad \lambda - \mu = b' \ . 
\end{displaymath}  
Along the first direction, $\ell(\theta)$ decreases rapidly to a local
unique minimum. The plot scale needs to be adapted in order to see
that a local unique minimum also exists along the second axis,
although the variations of $\ell(\theta)$ are much smaller along this
axis. It is therefore expected that any iterative minimization
algorithm will be directed rapidly to small values of $\ell(\theta)$
along the first axis. Reaching the global minimum by progressing along
the second axis will take much more steps. Moreover, it will also
require a precise evaluation of the increments of $\ell(\theta)$
between two neighboring values. We are then led to focus on the
variance of the values given by the Monte Carlo methods.
Figure~\ref{fig.likelihoodmc} shows box--plots of ~$200$ evaluations
of both methods at three different locations. Again we find that the
Brownian bridge method performs much better that the Pedersen method,
but both of them exhibit random fluctuations with an order of
magnitude noticeably greater than the variation of $\ell(\theta)$
along the second axis.

\subsection{Estimation results}
\begin{figure}
\begin{center}
\includegraphics[width=0.48\linewidth]{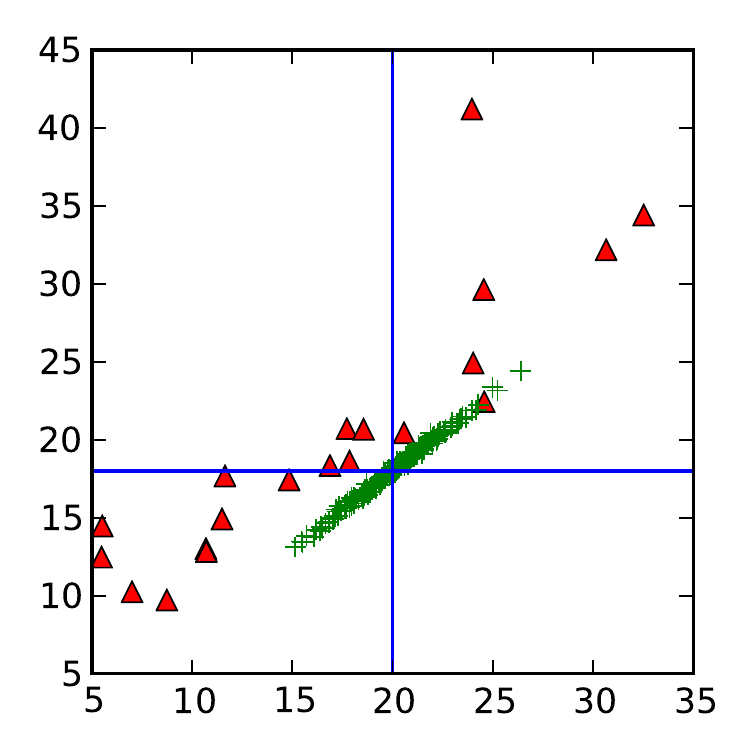} 
\includegraphics[width=0.48\linewidth]{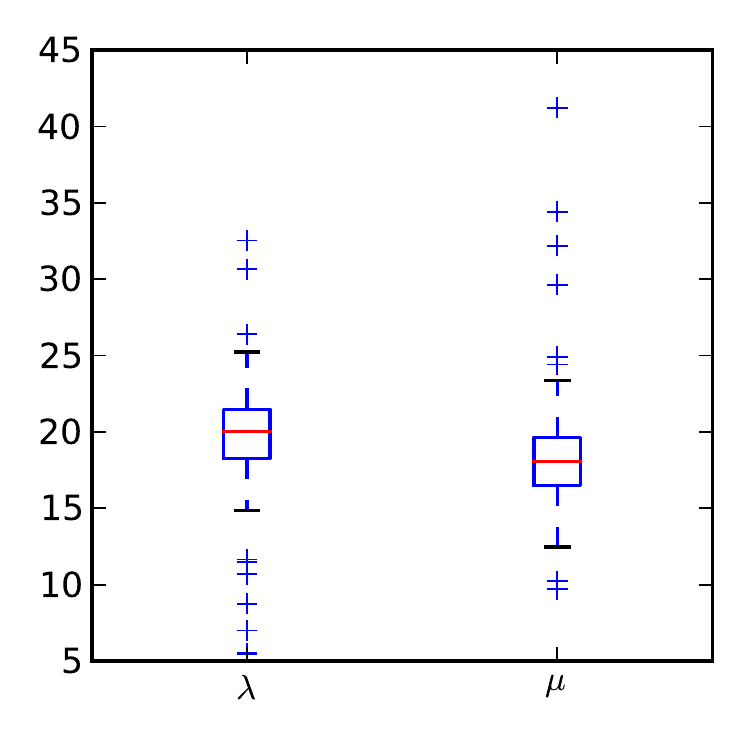}
\end{center} 
\caption{\it Empirical distribution of the MLE, using the finite
  difference approximation and parameters of Table~\ref{tab.par}. The
  right plot shows the marginal  distribution. On the left plot, the
  red triangles correspond to trajectories becoming extinct before
  final time, see Figure~\ref{fig.mc1}. Late extinction (the red
  triangle among green crosses) does not affect the estimation}  
\label{fig.est}
\end{figure}

%
As noticed in the previous section, the random fluctuations of the
Monte Carlo methods will mislay the iterative minimization
algorithm. The estimated value returned will be a local minimum,
depending essentially on the initial condition and the stopping
criterion. The estimation results of this section will then use the
finite difference approximation. On Figure~\ref{fig.est}, we plot the
estimated values for the data of Figure~\ref{fig.mc1}. We again notice
the asymmetry of the empirical distribution of the  MLE. The variance
is much smaller along the second axis, which means that information on
$(\lambda -\mu)$ in the data is better understood than information on
$(\lambda + \mu)$. In model~\eqref{eq.sde}, the drift coefficient
depends on $(\lambda -\mu)$ whereas the diffusion coefficient depends
on $(\lambda +\mu)$. In our case, inferring on the diffusion appears
to be the hardest task.

Finally, we observe the bad quality of the estimation when the data is
a  trajectory for which extinction happens during the transient
phase. Indeed, even if  we have taken good care of the extinction
problem, the information given by the data is not sufficient to allow
a useful inference on the parameter. On the other hand, the estimation is
correct if the observed trajectory reaches the stationary phase, even
if extinction occurs later on.   

\paragraph{Computational issues.}
The results presented herein have been obtained with a C++ code
using \verb+NLopt nonlinear-optimization package+
of Steven~G. Johnson, freely available at
\url{http://ab-initio.mit.edu/nlopt++}  which provides the
implementation of various algorithms. We have chosen to use a variant
of the Nelder and Mead simplex algorithm described in~\citet{rowan1990a}.

\paragraph{Acknowledgements.}
This work  was partially supported by the Laboratory of Excellence (Labex) NUMEV (Digital and Hardware Solutions, Modelling for the Environment and Life Sciences) coordinated by University of Montpellier 2, France.

\section{Concluding remarks}
%


The parameter estimation problem for a discretely observed diffusion
subject to almost sure extinction has been studied. In most practical
situations, the observed trajectory does not reach extinction. Indeed,
if $\lambda > \mu$, $\rho$ is small and $X_0$ is sufficiently large,
the mass absorbed at~$0$ is negligible, so that the transition kernel
is essentially $p_\Delta^\theta (y\,\vert\,x)\, \rmd y$. Nonetheless,
it makes sense to take care of the extinction probability
$E_\Delta(x)$ since the two parts are strongly related. Any
maximization algorithm will evaluate this likelihood for many values
of the parameter. For some of these values, neglecting the extinction
could lead to an abnormally high value of the likelihood and thus
mislay the maximization algorithm. Coupling extinction and
non--extinction in the complete Fokker--Planck equation results in a
more robust estimation procedure. This approach is not specific to the
particular application under consideration.  

Beyond the extinction problem and separate from it, the
identifiability of the parameters is also central in the model
presented here. With numerical experiments, we have seen that the
transport dynamics, given by the drift coefficient brings only a
partial information which does not allow to discriminate clearly
between the parameters values. The information on the demographic
noise encapsulated in the diffusion coefficient completes the
information, although it is much more difficult to extract from the
data. The Monte Carlo methods presented here are not able to achieve
the precision level required to that end.

\clearpage
\appendix
\section*{Appendix}
\section{Derivation of the stochastic logistic model \eqref{eq.sde}}
\label{appendix.sde}

A natural way to derive the SDE  \eqref{eq.sde} is to consider, at a microscopic scale, a population size $N_{t}$ subject to a birth and death process that features a logistic mechanism, e.g.:
\begin{align}
\label{eq.proc.N_t}
  \P(N_{t+h}=n'|N_{t}=n)
  \equiva_{h\to 0}
  \begin{cases}
    h\, \lambda\,n + o(h)\,, & \textrm{if } n'=n+1\,,
  \\
    \textstyle
    h\, (\mu+\frac{\alpha}{\kappa}\,n)\,n+o(h)\,, & \textrm{if } n'=n-1\,,
  \\
    \textstyle
    1-h\, (\lambda+\mu+\frac{\alpha}{\kappa}\,n)\,n+o(h)\,, & \textrm{if } n'=n\,,\\
    o(h)\,, & \textrm{otherwise}
  \end{cases}
\end{align}
with $\lambda>0$, $\mu>0$, $\alpha>0$, $\kappa>0$; here the birth per capita rate $\lambda$ is constant and the death per capita rate $\mu+\frac{\alpha}{\kappa}\,n$ increases linearly with the population size $n$.

The rescaled process:
\[
  X^\kappa_{t}
  \eqdef
  \frac1\kappa\,N_{t}
\]
is a pure jump Markov process with values in $\frac1\kappa\,\N$ where $\kappa$
denotes the order of magnitude of the population size $N_{t}$ of interest.
The rescaled process is a pure jump Markov process and its distribution law 
is characterized by its infinitesimal generator:
\begin{align*} 
\label{eq.gen:N}
 \lim_{t\to 0}
 \textstyle
 \frac{\E(\phi(X^\kappa_{t})|X^\kappa_{0}=x)-\phi(x)}{t} 
  = 
   \lambda \, x \, \left [
     \phi (x + \frac{1}{\kappa}) - \phi (x)
     \right ] 
   + (\mu + \alpha \, x) \, x \, \left [
     \phi (x - \frac{1}{\kappa}) - \phi (x)
   \right ] 
\end{align*} 
defined for any bounded function $\phi$.
In large population scale, that is for $\kappa$ large, 
a {\em diffusion approximation} of
$X^\kappa_t$ is obtained by performing a second order Taylor expansion of
regular functions $\phi$. This yields the operator $\AAA$ defined by
\eqref{eq.gen:X}, which is the infinitesimal generator of the
diffusion process solution of~\eqref{eq.sde}.

\medskip

Note
that~\eqref{eq.sde} can be rewritten as
\begin{displaymath}
  \rmd X_t = r\, \left ( 1 - \frac{X_t}{K} \right ) \, X_t\, \rmd t
  + \rho \, \sqrt{ r'\,  \left ( 1 + \frac{X_t}{K'} \right ) \, X_t } 
  \,\rmd B_t
\end{displaymath} 
with $r = \lambda - \mu,$ $K = \frac{r}{\alpha}$, $r' = \lambda +
\mu$, $K' = \frac{r'}{\alpha}$ and $\rho=\frac1{\sqrt{\kappa}}$, so that
its instantaneous mean has the 
same form as in~\eqref{eq.verhulst}. However, the deterministic
model~\eqref{eq.verhulst} does not describe the evolution of
$\E[X_t]$. Diffusion approximation technique is frequently encountered
in life sciences. A typical example of a two--dimensional bioreactor
is described in~\citet{campillo2011chemostat} or~\citet{joannides2013a}.

\medskip

The relationship between  the distribution laws of the process
 $X^\kappa_{t}$ and its \emph{diffusion approximation} $X_{t}$ for large
$\kappa$ are precised in~\citet[Chap.~7
  and~11]{ethier1986a}.

\section{Existence and uniqueness for SDE \eqref{eq.sde}}
\label{appendix.sde.existence.uniqueness}

The drift function $b$ is locally Lipschitz on $\R$ but
fails to satisfy the usual linear growth condition.  On the other
hand, the diffusion function $\sigma$ is not locally Lipschitz on
$\R$. Nevertheless, we have
\begin{lemma}
For any non-negative initial condition $X_0 \in L^2$, there exists a
unique non-negative solution to~\eqref{eq.sde}. 
\end{lemma} 

\proof
We first deal with the diffusion term. Consider $h:\R \mapsto\R$
globally Lipschitz and of linear growth, suppose also that $h(0) =
0$. Then introduce for $\ell \geq 1$
\begin{displaymath} 
\sigma _\ell (x) \eqdef
\begin{dcases}
  \sigma (x)\,,
  & \text{if }  \textstyle \frac{1}{\ell}\leq x\,,
\\[0.5em]
  \textstyle \rho\, \sqrt{\frac{1}{\ell}\, (\lambda + \mu + 
                \frac{\alpha}{\ell}x)\, (2\, \ell \, x - 1)}  \,,
  & \text{if } \textstyle \frac{1}{2\,\ell} < x < \frac{1}{\ell}\,,
\\[0.5em] 
  0\,, 
  &\textstyle\text{if } x\leq \frac{1}{2\,\ell}\,.
\end{dcases} 
\end{displaymath} 
This function is globally Lipschitz so that SDE
\begin{displaymath}
  \rmd Y_t^\ell  
  = 
  h(Y_t^\ell)\, \rmd t + \sigma_\ell (Y_t^\ell)\, \rmd B_t
\end{displaymath} 
has a unique solution with a.s.\! continuous path, for any non-negative
initial condition $Y_0 \in L^2$. Define $T_\ell = \inf\{t
\geq 0\,; \,  Y_t^\ell \leq \frac{1}{\ell} \}$ and note that $Y_t^\ell$ and
$Y_t^{\ell'}$ coincide up to time $T_{\ell'}$ when ${\ell'} \leq \ell$~(Durett, Lemma 2.8,
Chap.5). The process $Y_t = Y_t^\ell$ is then well defined up to time
$T_\infty = \lim_{\ell\uparrow \infty} T_{\ell}$. Now since $\sigma _\ell$   and
$\sigma$ coincide on $[\frac{1}{\ell}, +\infty)$, $Y_t$ is a solution to
the following SDE 
\begin{equation}
\label{eq.pb.in.0} 
\rmd Y_t  = h(Y_t)\, \rmd t + \sigma (Y_t)\, \rmd B_t
\end{equation} 
on time interval $[0,T_\infty)$. 
On the event $\{T_\infty < \infty\}$, we define $Y_t = 0$ for $t\geq
T_\infty$, which is an obvious solution of~\eqref{eq.pb.in.0} with
initial condition~$Y_0 = 0$. We conclude that~\eqref{eq.pb.in.0} has a
unique non-negative solution for all $t\geq 0$ and non-negative initial
condition $Y_0\in L^2$. 

Likewise, let for $\ell \geq 1$
\begin{displaymath} 
b _\ell (x) \eqdef
\begin{dcases}
  b (x)\,,
  &\text{if }  x \leq \ell \,,
\\
  (\lambda - \mu - \alpha \, \ell\,x)\, (2\, \ell - x)\,,
  &
  \text{if } \ell < x < 2\, \ell \,,
\\
  0 \,,
  &
  \text{if }2\,\ell\leq x\,.
\end{dcases} 
\end{displaymath} 
Since $b_\ell$ is globally Lipschitz and bounded, the preceding result applies so
that
\begin{displaymath} 
\rmd X_t^\ell = b_\ell(X_t^\ell)\, \rmd t + \sigma(X_t^\ell)\,\rmd B_t
\end{displaymath} 
has a unique non-negative solution for each $\ell$. Consider the stopping
times $S_\ell = \inf \{ t\geq 0\,; \, X_t^\ell \geq \ell\}$ and 
$S_\infty = \lim_{\ell\uparrow \infty} S_\ell$, and define $X_t =X_t^\ell$ for
$t \leq S_\infty $. Since $b_\ell$ and $b$ coincide on $[0,\ell]$, it holds
\begin{displaymath}
X_{t\wedge S_\ell}  = X_0 + \int_0^{t \wedge S_\ell} b(X_s)\, \rmd s + 
 \int_0^{t \wedge S_\ell} \sigma (X_s)\, \rmd B_s \,.
\end{displaymath} 
The last term is a martingale, so by the
optional stopping theorem
\begin{align*} 
\E(X_{t\wedge S_\ell}) 
&= 
  \E(X_0) 
  + 
  \E\left (\int_0^{t \wedge S_\ell} b(X_s)\,\rmd s \right) 
\\
&= 
  \E(X_0) 
  + 
  \int_0^t \E [\indic_{(s\leq S_\ell)}\, b(X_{s\wedge S_\ell})] \,\rmd s 
\\
& \leq 
  \E(X_0) 
  + \int_0^t \E[ b(X_{s\wedge S_\ell})] \,\rmd s 
\\
& \leq  
  \E(X_0) 
  + \int_0^t \E[ (\lambda - \mu)\, X_{s\wedge S_\ell} ]\, \rmd s 
  \leq 
  \E(X_0) + (\lambda - \mu)\, \int_0^t \E[ X_{s\wedge S_\ell} ] \,\rmd s \,.
\end{align*} 
The Gronwall lemma  yields $\E(X_{T\wedge S_\ell}) \leq \E(X_0) \,
\exp\{ (\lambda - \mu )\, T\}$, for all $T >0$, and from Fatou's lemma:
\begin{displaymath} 
\forall \, T > 0,\qquad \E(X_{T\wedge S_\infty}) \leq \liminf_{\ell\to\infty}
\E(X_{T\wedge S_\ell}) \leq \exp\{ (\lambda - \mu )\, T\}
\end{displaymath} 
so that $\P(S_\infty \leq T) = 0$. Hence $S_\infty = \infty$ p.s.\! and
the lemma is proved.
\hfill $\square$

\section{Algorithms}

\paragraph{Notations}
\begin{itemize}
\item $\mathrm D(m,s^2,x)$: probability density function of $\NN(m,s^2)$
  at point $x$
\item  $\mathrm P(m,s^2,p)$: cumulative distribution function $\P(\NN(m,s
  ^2)\leq p)$ at point $p$
\item  $\mathrm R(m,s^2)$: random deviate of $\NN(m,s ^2)$
\end{itemize}

\begin{algorithm}
\caption{This algorithm returns a random value for $\qP_\Delta
  (y\,\vert\,x)$, the approximation obtained with the  Pedersen
  method as described in Section~\ref{sec.pedersen}}
\label{algo:ped}
\begin{algorithmic}
\REQUIRE Current observation $x>0$
\REQUIRE Next observation $y>0$
\FOR{$i = 1 \to N$} 
\STATE $\xi \gets x$ 
\STATE $t \gets 0$
    \FOR{$j = 1 \to n - 1$} 
    \IF{$\xi \not = 0$} 
        \STATE $\mathrm{mean} \gets \xi + \delta\, b(\xi)$
        \STATE $\mathrm{sd} \gets \sigma(\xi) $
        \STATE $\xi \gets \max \left \{0, \mathrm{mean} + \sqrt{\delta} \,
       \mathrm{sd} * \mathrm{R}(0,1)\right \}$
    \ENDIF
    \STATE $t \gets t + \delta$
    \ENDFOR
    \IF{$y \not = 0$} 
    \STATE  $v \gets \mathrm{D}(\xi + \delta \,
    b(\xi),\delta \sigma ^2(\xi),y)$
    \ELSIF{$\xi  = 0$} 
    \STATE $v \gets 1$
    \ELSE
    \STATE $v \gets \mathrm{P}(\xi + \delta \,
    b(\xi),\delta \sigma ^2(\xi),0)$
    \ENDIF
    \RETURN $\displaystyle \frac{v}{N}$.
\ENDFOR

\end{algorithmic}
\end{algorithm}

\begin{algorithm}
\caption{This algorithm returns a random value for $\qB_\Delta
  (y\,\vert\,x)$, the approximation obtained with the Brownian bridge
  variant described in text, Section~\ref{sec.bb}.}
\label{algo:bb}
\begin{algorithmic}
\REQUIRE Current observation $x>0$
\REQUIRE Next observation $y>0$
\STATE $L \gets 0$
\FOR{$i = 1 \to N$} 
\STATE $\xi \gets x$ 
\STATE $t \gets 0$
\STATE $\psi \gets 1$
    \FOR{$j = 1 \to n - 1$} 
    \IF{$\xi \not = 0$} 
        \STATE $\mathrm{euler.mean} \gets \xi + \delta\, b(\xi)$
        \STATE $\mathrm{bridge.mean} \gets \xi + \delta\,
        \displaystyle \frac{y -
          \xi}{\Delta - t}$
        \STATE $\mathrm{sd} \gets \sigma(\xi) $
        \STATE $\xi \gets \max \left \{0, \mathrm{bridge.mean} + \sqrt{\delta} \,
       \mathrm{sd} * \mathrm{R}(0,1)\right \}$
        \IF{$\xi \not = 0$} 
        \STATE $\psi \gets \psi \, 
        \displaystyle \frac{\mathrm{D}(\mathrm{euler.mean},\mathrm{sd}^2,\xi)}
             {\mathrm{D}(\mathrm{bridge.mean},\mathrm{sd}^2,\xi)}$
        \ELSE
        \STATE $\psi \gets \psi \, \displaystyle \frac{
          \mathrm{P}(\mathrm{euler.mean},\mathrm{sd}^2,0)}
               {\mathrm{P}(\mathrm{bridge.mean},\mathrm{sd}^2,0)}$
        \ENDIF
    \ENDIF
    \STATE $t \gets t + \delta$
    \ENDFOR
    \IF{$y \not = 0$} 
    \STATE $L \gets L + \psi * \mathrm{D}(\xi + \delta \,
    b(\xi),\delta \sigma ^2(\xi),y)$
    \ELSIF{$\xi  = 0$} 
    \STATE  $L \gets L + \psi$
    \ELSE
    \STATE $L \gets L + \psi * \mathrm{P}(\xi + \delta \,
    b(\xi),\delta \sigma ^2(\xi),y)$
    \ENDIF
\ENDFOR
\RETURN $\displaystyle \frac{L}{N}$

\end{algorithmic}
\end{algorithm}


\cleardoublepage
\addcontentsline{toc}{section}{Reference}
\bibliographystyle{plainnat}

\bibliography{hal}

\end{document}

%% file: discret.pdf_tex
\begingroup%
  \makeatletter%
  \providecommand\color[2][]{%
    \errmessage{(Inkscape) Color is used for the text in Inkscape, but the package 'color.sty' is not loaded}%
    \renewcommand\color[2][]{}%
  }%
  \providecommand\transparent[1]{%
    \errmessage{(Inkscape) Transparency is used (non-zero) for the text in Inkscape, but the package 'transparent.sty' is not loaded}%
    \renewcommand\transparent[1]{}%
  }%
  \providecommand\rotatebox[2]{#2}%
  \ifx\svgwidth\undefined%
    \setlength{\unitlength}{523.11640625bp}%
    \ifx\svgscale\undefined%
      \relax%
    \else%
      \setlength{\unitlength}{\unitlength * \real{\svgscale}}%
    \fi%
  \else%
    \setlength{\unitlength}{\svgwidth}%
  \fi%
  \global\let\svgwidth\undefined%
  \global\let\svgscale\undefined%
  \makeatother%
  \begin{picture}(1,0.27594539)%
    \put(0,0){\includegraphics[width=\unitlength]{discret.pdf}}%
    \put(0.1131739,0.25503107){\color[rgb]{0,0,1}\makebox(0,0)[lb]{\smash{$\Delta$}}}%
    \put(0.02141612,0.00269717){\color[rgb]{1,0,0}\makebox(0,0)[lb]{\smash{$\delta$}}}%
    \put(-0.00152332,0.17091977){\color[rgb]{0,0,1}\makebox(0,0)[lb]{\smash{$\tto_0=0$}}}%
    \put(-0.00152332,0.07151551){\color[rgb]{1,0,0}\makebox(0,0)[lb]{\smash{$t_0=0$}}}%
    \put(0.24316409,0.07151551){\color[rgb]{1,0,0}\makebox(0,0)[lb]{\smash{$t_n=\Delta$}}}%
    \put(0.15140631,0.07151551){\color[rgb]{1,0,0}\makebox(0,0)[lb]{\smash{$t_k$}}}%
    \put(0.24316409,0.17091977){\color[rgb]{0,0,1}\makebox(0,0)[lb]{\smash{$\tto_1$}}}%
    \put(0.88546855,0.07151551){\color[rgb]{0,0,0}\makebox(0,0)[lb]{\smash{time $t$}}}%
    \put(0.73253892,0.17091977){\color[rgb]{0,0,1}\makebox(0,0)[lb]{\smash{$\tto_M=T$}}}%
    \put(0.73253892,0.07151551){\color[rgb]{0,0,0}\makebox(0,0)[lb]{\smash{$T$}}}%
  \end{picture}%
\endgroup%

%% file: pedersen.pdf_tex
\begingroup%
  \makeatletter%
  \providecommand\color[2][]{%
    \errmessage{(Inkscape) Color is used for the text in Inkscape, but the package 'color.sty' is not loaded}%
    \renewcommand\color[2][]{}%
  }%
  \providecommand\transparent[1]{%
    \errmessage{(Inkscape) Transparency is used (non-zero) for the text in Inkscape, but the package 'transparent.sty' is not loaded}%
    \renewcommand\transparent[1]{}%
  }%
  \providecommand\rotatebox[2]{#2}%
  \ifx\svgwidth\undefined%
    \setlength{\unitlength}{612.02680664bp}%
    \ifx\svgscale\undefined%
      \relax%
    \else%
      \setlength{\unitlength}{\unitlength * \real{\svgscale}}%
    \fi%
  \else%
    \setlength{\unitlength}{\svgwidth}%
  \fi%
  \global\let\svgwidth\undefined%
  \global\let\svgscale\undefined%
  \makeatother%
  \begin{picture}(1,0.48033422)%
    \put(0,0){\includegraphics[width=\unitlength]{pedersen.pdf}}%
    \put(0.64956741,0.27466748){\color[rgb]{0,0,1}\makebox(0,0)[lb]{\smash{\scriptsize$\bar X^{(i)}_{\Delta-\delta}$}}}%
    \put(0.31382671,0.00345803){\color[rgb]{0,0,0}\makebox(0,0)[lb]{\smash{$0$}}}%
    \put(0.72452711,0.15275617){\color[rgb]{0,0,0}\makebox(0,0)[lb]{\smash{$\Delta-\delta$}}}%
    \put(0.81496233,0.18563122){\color[rgb]{0,0,0}\makebox(0,0)[lb]{\smash{$\Delta$}}}%
    \put(0.19120084,0.10131368){\color[rgb]{0,0,0}\makebox(0,0)[lb]{\smash{$x$}}}%
    \put(0.72817135,0.37006005){\color[rgb]{0,0,1}\makebox(0,0)[lb]{\smash{$K_\delta(\bar X^{(i)}_{\Delta-\delta},\rmd y)$}}}%
    \put(0.70276827,0.43809384){\color[rgb]{0,0,0}\makebox(0,0)[lb]{\smash{$\QP_\Delta(x,\rmd y)$}}}%
  \end{picture}%
\endgroup%

%% file: bridge.pdf_tex
\begingroup%
  \makeatletter%
  \providecommand\color[2][]{%
    \errmessage{(Inkscape) Color is used for the text in Inkscape, but the package 'color.sty' is not loaded}%
    \renewcommand\color[2][]{}%
  }%
  \providecommand\transparent[1]{%
    \errmessage{(Inkscape) Transparency is used (non-zero) for the text in Inkscape, but the package 'transparent.sty' is not loaded}%
    \renewcommand\transparent[1]{}%
  }%
  \providecommand\rotatebox[2]{#2}%
  \ifx\svgwidth\undefined%
    \setlength{\unitlength}{581.53071289bp}%
    \ifx\svgscale\undefined%
      \relax%
    \else%
      \setlength{\unitlength}{\unitlength * \real{\svgscale}}%
    \fi%
  \else%
    \setlength{\unitlength}{\svgwidth}%
  \fi%
  \global\let\svgwidth\undefined%
  \global\let\svgscale\undefined%
  \makeatother%
  \begin{picture}(1,0.5175704)%
    \put(0,0){\includegraphics[width=\unitlength]{bridge.pdf}}%
    \put(0.66109543,0.30509752){\color[rgb]{0,0,1}\makebox(0,0)[lb]{\smash{\scriptsize$\tilde X^{(i)}_{\Delta-\delta}$}}}%
    \put(0.33028412,0.00363936){\color[rgb]{0,0,0}\makebox(0,0)[lb]{\smash{$0$}}}%
    \put(0.76252209,0.16076685){\color[rgb]{0,0,0}\makebox(0,0)[lb]{\smash{$\Delta-\delta$}}}%
    \put(0.85769983,0.1953659){\color[rgb]{0,0,0}\makebox(0,0)[lb]{\smash{$\Delta$}}}%
    \put(0.20122762,0.10662667){\color[rgb]{0,0,0}\makebox(0,0)[lb]{\smash{$x$}}}%
    \put(0.81124033,0.36428814){\color[rgb]{0,0,1}\makebox(0,0)[lb]{\smash{$K_\delta(\tilde X^{(i)}_{\Delta-\delta},\rmd y)$}}}%
    \put(0.7589697,0.46083963){\color[rgb]{0,0,0}\makebox(0,0)[lb]{\smash{$\QB_\Delta(x,\rmd y)$}}}%
    \put(0.66928855,0.38963066){\color[rgb]{0,0,0}\makebox(0,0)[lb]{\smash{\scriptsize $y$}}}%
    \put(0.56961966,0.22351557){\color[rgb]{0,0,1}\makebox(0,0)[lb]{\smash{\scriptsize trajectory}}}%
    \put(0.56856699,0.18110582){\color[rgb]{0,0,1}\makebox(0,0)[lb]{\smash{\scriptsize $\tilde X_{t_{0:n-1}}^{(i)}$}}}%
  \end{picture}%
\endgroup%